\documentclass[11pt]{article}
\usepackage[final]{epsfig}
\usepackage{graphics}
\usepackage{amsmath}
\usepackage{amsfonts}
\usepackage{latexsym}
\usepackage{amssymb}
\usepackage{graphicx}
\usepackage{epstopdf}
\usepackage{amsthm,amscd,array,stmaryrd,mathrsfs}

\let\ssection=\section
\renewcommand{\section}{\setcounter{equation}{0}\ssection}

\theoremstyle{plain}
\newtheorem{theorem}{Theorem}
\newtheorem{lemma}{Lemma}[section]
\newtheorem{corollary}[lemma]{Corollary}

\newtheorem{proposition}[lemma]{Proposition}
\newtheorem{remark}[lemma]{Remark}
\newtheorem{example}[lemma]{Example}

\begin{document}
\newcommand{\eps}{{\varepsilon}}
\def\a{\alpha}
\def\b{\beta}
\def\t{\tau}
\newcommand{\g}{{\gamma}}
\newcommand{\G}{{\Gamma}}
\newcommand{\proofend}{$\Box$\bigskip}
\newcommand{\C}{{\mathbb C}}
\newcommand{\cC}{\mathcal{C}}
\newcommand{\Q}{{\mathbb Q}}
\newcommand{\R}{{\mathbb R}}
\newcommand{\Z}{{\mathbb Z}}
\newcommand{\N}{{\mathbb N}}
\newcommand{\RP}{{\mathbb{RP}}}
\newcommand{\CP}{{\mathbb{CP}}}
\newcommand{\pP}{{\mathbb{P}}}
\newcommand{\cP}{{\mathcal{P}}}
\newcommand{\PSL}{{\mathrm{PSL}}}
\newcommand{\PGL}{{\mathrm{PGL}}}
\newcommand{\SL}{{\mathrm{SL}}}
\def\proof{\paragraph{Proof.}}
\def\startproof{\paragraph{Proof.}}

\newcounter{vo}
\newcommand{\vo}[1]
{\stepcounter{vo}$^{\bf VO\thevo}$%
\footnotetext{\hspace{-3.7mm}$^{\blacksquare\!\blacksquare}$
{\bf VO\thevo:~}#1}}

\newcounter{st}
\newcommand{\st}[1]
{\stepcounter{st}$^{\bf ST\thest}$%
\footnotetext{\hspace{-3.7mm}$^{\blacksquare\!\blacksquare}$
{\bf ST\thest:~}#1}}

\newcounter{rs}
\newcommand{\rs}[1]
{\stepcounter{rs}$^{\bf RS\thers}$%
\footnotetext{\hspace{-3.7mm}$^{\blacksquare\!\blacksquare}$
{\bf RS\thers:~}#1}}

\title{Liouville-Arnold integrability of the pentagram map on closed polygons}

\author{Valentin Ovsienko
\and
Richard Evan Schwartz
\and
Serge Tabachnikov}

\date{}
\maketitle

\begin{abstract}
The pentagram map is a discrete dynamical system defined on the moduli space
of polygons in the projective plane.
This map has recently attracted a considerable interest, mostly because
its connection to a number of different domains, such as:
classical projective geometry, algebraic combinatorics,
moduli spaces, cluster algebras and integrable systems.

Integrability of  the pentagram map was conjectured in \cite{Sch2}
and proved in \cite{OST} for a larger space of twisted polygons.
In this paper, we prove the initial conjecture
that the pentagram map is completely integrable
on the moduli space of closed polygons.
In the case of convex polygons in the real projective plane, this result
implies the existence of a toric foliation on the moduli space.
The leaves of the foliation carry affine structure and the dynamics of the pentagram map
is quasi-periodic.
Our proof is based on an invariant Poisson structure on the space of twisted polygons.
We prove that the Hamiltonian vector fields corresponding to the monodoromy invariants
preserve the space of closed polygons and define an invariant affine structure on the level surfaces of the monodromy invariants.
\end{abstract}

\thispagestyle{empty}

\tableofcontents

\section{Introduction}
\begin{figure}[hbtp]
\centering
\includegraphics[height=150pt]{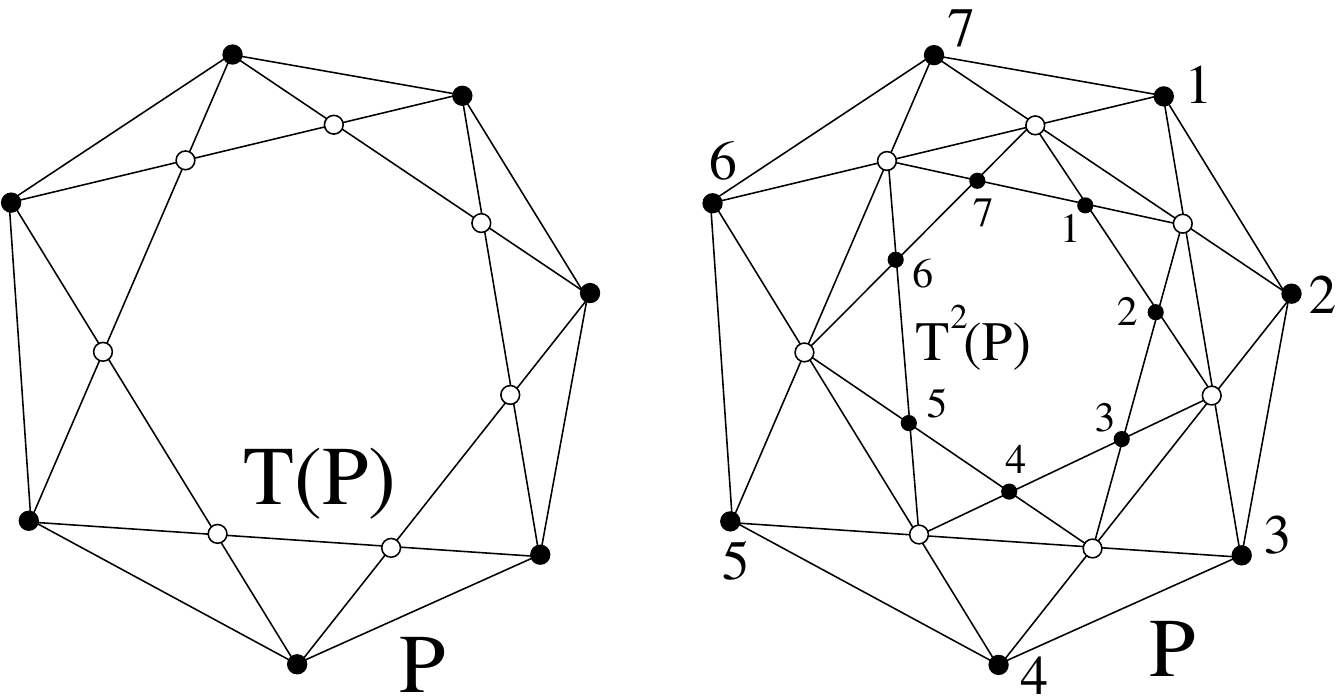}
\newline
\caption{The pentagram map and its second iterate defined on a convex $7$-gon}
\label{5and6}
\end{figure}

The {\it pentagram map} is a geometric 
construction which carries one polygon to another.
Given an $n$-gon $P$, the vertices of
the image $T(P)$ under the pentagram map
are the intersection points of consecutive shortest diagonals of~$P$.
The left side of Figure \ref{5and6} shows the basic construction.
The right hand side shows the second iterate of the 
pentagram map. The second iterate has the virtue that it acts
in a canonical way on a labeled polygon, as indicated.
The first iterate also acts on labeled polygons, but one must
make a choice of labeling scheme; see Section \ref{LEFT}.
The simplest example of the pentagram map for pentagons was considered in \cite{Mot}.
In the case of arbitrary $n$, the map was introduced in \cite{Sch1} and 
further studied in \cite{Sch2,Sch3}.

The pentagram map is defined on any polygon whose points are
in general position, and also on some polygons whose points are not in general
position.  One sufficient condition for the pentagram map to be
well defined is that every consecutive triple of points is not collinear.
However, this last condition is not invariant under the pentagram map.

The pentagram map commutes with projective transformations
and thereby induces a (generically defined) map
\begin{equation}
\label{TheMap}
T:\cC_n\to\cC_n
\end{equation}
where $\cC_n$ is the moduli space of projective
equivalence classes of $n$-gons in the projective plane.
Mainly we are interested in the subspace $\cC_n^0$ of
projective classes convex $n$-gons.  The
pentagram map is entirely defined on
$\cC_n^0$ and preserves this subspace.

Note that the pentagram map can be defined over an arbitrary field.
Usually, we restrict our considerations to the geometrically natural real case
of convex $n$-gons in $\RP^2$.
However,  the complex case represents a special interest since
the moduli space of $n$-gons in $\CP^2$ is a higher analog of the moduli space
$\mathcal{M}_{0,n}$.
Unless specified,
we will be using the general notation $\pP^2$ for the projective plane
and $\PGL_3$ for the group of projective transformations.

\subsection{Integrability problem and known results}

\begin{figure}[hbtp]
\centering
\includegraphics[height=150pt]{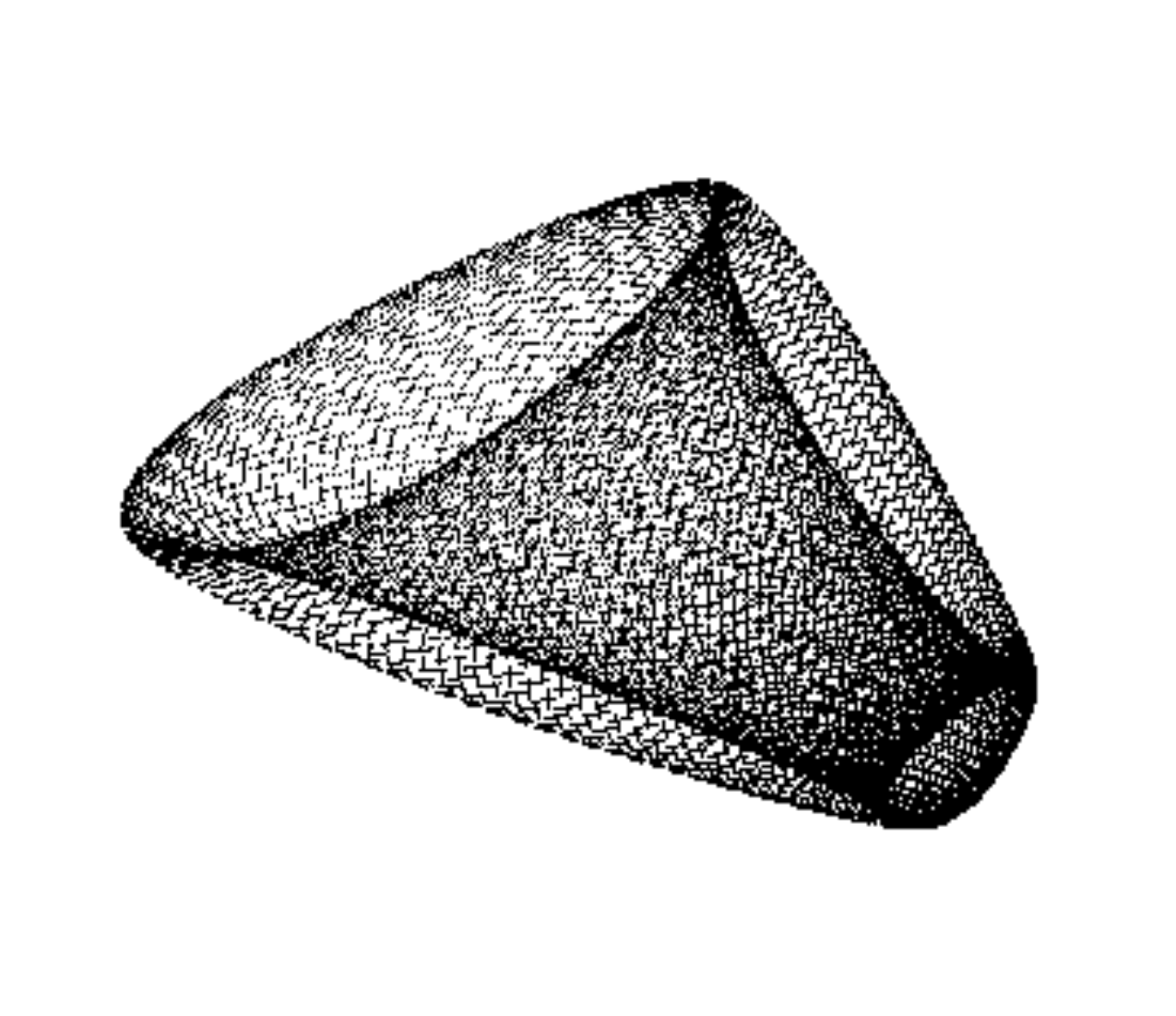}
\newline
\caption{An orbit of the pentagram map on a heptagon}
\label{Torus}
\end{figure}

Assuming that the labeling schemes have been chosen
carefully, the map $T:\cC_5\to\cC_5$ is the identity map and
the map $T:\cC_6\to\cC_6$ is  an involution.  
See \cite{Sch1}.
The conjecture that the map (\ref{TheMap}) is completely integrable
was formulated roughly in \cite{Sch1} and then
more precisely in \cite{Sch2}.  This conjecture was inspired
by computer experiments in the case $n=7$.
Figure \ref{Torus} presents (a two-dimensional projection of) an orbit of a convex heptagon in $\RP^2$.

The first results regarding the integrability of the pentagram map were proved for
the pentagram map defined on a larger space,~$\cP_n$, of \textit{twisted} $n$-gons.
A series of $T$-invariant functions (or first integrals)
called the \textit{monodromy invariants},
was constructed in \cite{Sch3}.
In \cite{OST} (see also \cite{OST1} for a short version),
the complete integrability of $T$ on $\cP_n$ was proved
with the help of a $T$-invariant Poisson structure,
such that the monodromy invariants Poisson-commute.

In~\cite{Sol}, F. Soloviev found a Lax representation of  the pentagram map
and proved its algebraic-geometric integrability. 
The space of polygons (either $\cP_n$ or $\cC_n$) 
is parametrized in terms of a spectral curve with marked points and a divisor. 
The spectral curve is determined by the monodromy invariants, 
and the divisor corresponds to a point on a torus -- the Jacobi variety of the spectral curve.  These results allow one to construct explicit solutions formulas using Riemann theta functions  (i.e., the variables that determine the polygon as explicit functions of time).  Soloviev also deduces the invariant Poisson bracket of~\cite{OST}
from the Krichever-Phong universal formula.

Our result below has the same dynamical implications as that of Soloviev,
in the case of real convex polygons. Soloviev's approach is by way of 
algebraic integrability, and it has the advantage that it identifies the invariant tori explicitly as certain Jacobi varieties. Our proof is in the framework of Liuoville-Arnold integrability, and  it is more direct and self-contained.

\subsection{The main theorem}

The main result of the present paper is to give a purely geometric proof of the following
result.

\begin{theorem}
\label{Main}
Almost every point of $\cC_n$ lies on a $T$-invariant algebraic submanifold
of dimension 
\begin{equation}
\label{Dim}
d=\left\{
\begin{array}{l}n-4, \;n \; \hbox{is odd}\\
n-5, \; n \; \hbox{is even}.
\end{array}\right.
\end{equation}
 that has
a $T$-invariant affine structure.
\end{theorem}

Recall that an affine structure on a $d$-dimensional manifold
is defined by a locally free action of the $d$-dimensional Abelian Lie algebra,
that is, by $d$ commuting vector fields linearly independent at every point.

In the case of convex $n$-gons in the real projective plane,
thanks to the compactness of the space established in \cite{Sch2}, our result reads:

\begin{corollary}
\label{ac}
Almost every orbit in $\cC^0_n$ lies on a finite union
of smooth $d$-dimensional tori, where $d$ is
as in equation (\ref{Dim}).  The union of
these tori has 
a $T$-invariant affine structure.
\end{corollary}

\noindent
Hence, the orbit of almost every convex
$n$-gon undergoes quasi-periodic motion under the pentagram map.
The above statement is precisely the integrability theorem in the
Liouville--Arnold sense~\cite{Arn}.

Let us also mention that the dimension of the invariant sets given by (\ref{Dim})
is precisely a half of the dimension of $\cC_n$, provided $n$ is odd,
which is a usual, generic, situation for an integrable system.
If $n$ is even, then $d=\frac{1}{2}\dim\cC_n-1$ so that one can talk of
``hyper-integrability''.

Our approach is based on the results of \cite{Sch3} and \cite{OST}.
We prove that the level sets of the monodromy invariants on the subspace $\cC_n\subset\cP_n$
are algebraic subvarieties of $\cC_n$ of dimension~(\ref{Dim}).
We then prove that the Hamiltonian vector fields corresponding to the
invariant functions are tangent to $\cC_n$
(and therefore to the level sets).
Finally, we prove that the Hamiltonian vector fields define an 
affine structure on a generic level set.  The main calculation, which
establishes the needed independence of the monodromy invariants
and their Hamiltonian vector fields, uses a 
trick that is similar in spirit to tropical algebra.

One point that is worth emphasizing is that our proof does not
actually produce a symplectic (or Poisson) structure on the
space ${\cal C\/}_n$.  Rather, we use the Poisson structure
on the ambient space $\cP^n$, together with the invariants, to produce
enough commuting flows on ${\cal C\/}_n$ in order to fill
out the level sets.  

\subsection{Related topics}\label{Relation}

The pentagram map is a particular example of a discrete integrable system.
The main motivation for studying this map is its relations to
different subjects, such as: 
a) projective differential geometry;
b) classical integrable systems and symplectic geometry;
c) cluster algebras;
d) algebraic combinatorics of Coxeter frieze patterns.
All these relations may be beneficial not only for the study of the pentagram map,
but also for the above mentioned subjects.
Let us mention here some recent developments involving the pentagram map.

\begin{itemize}
\item The relation of~$T$ to the classical Boussinesq equation was essential for \cite{OST}.
In particular, the Poisson bracket was obtained as a discretization of the
(first) Adler-Gelfand-Dickey bracket related to the Boussinesq equation.
We refer to \cite{TN,LN} and references therein for more information about
different versions of the discrete Boussinesq equation.
\item In~\cite{ST1}, surprising results of elementary  projective geometry are obtained
in terms of the pentagram map, its iterations and generalizations.
\item In \cite{ST2}, special relations amongst the monodromy
invariants are established for polygons that are inscribed into
a conic. 
\item In \cite{FM}, the pentagram map is related to Lie-Poisson loop groups.
\item The paper \cite{MB} concerns discretizations of Adler-Gelfand-Dickey flows as multi-dimensional generalizations of the pentagram map.
\item A particularly interesting feature of the pentagram map is its relation to
the theory of cluster algebras developed by Fomin and Zelevinsky, see \cite{FZ1}.
This relation was noticed in \cite{OST} and developed
in \cite{Gli}, where the pentagram map on the space of twisted $n$-gons is interpreted as a
sequence of cluster algebra mutations, and an explicit formula for the iterations of $T$
is calculated\footnote{This can be understood as a version of integrability
or ``complete solvability''.}.
\item The structure of cluster manifold on the space $\cC_n$ and the related notion of $2$-frieze pattern are
investigated in \cite{SVS}.
\end{itemize}

\section{Integrability on the space of twisted $n$-gons}\label{TwiS}

In this section, we explain the proof of the main result in
our paper \cite{OST},
the Liouville-Arnold integrability of the pentagram map on the space of twisted $n$-gons.
While we omit some technical details, we take the opportunity to fill a gap
in \cite{OST}:  there we claimed that the monodromy invariants
Poisson commute, but our proof there had a flaw.  Here we present a correct
 proof of this fact.

\subsection{The space $\cP_n$}

We recall the definition of the space of twisted $n$-gons.

A \textit{twisted $n$-gon} is a
map $\phi: \Z \to\pP^2$ such that
\begin{equation} 
\label{PoT}
v_{i+n}=M \circ v_i,
\end{equation}
for all $i \in \Z$ and some fixed element $M\in\PGL_3$ called the \textit{monodromy}.
We denote by $\cP_n$
the space of twisted $n$-gons modulo
projective equivalence.
The pentagram map extends to a generically
defined map $T:\cP_n\to\cP_n$.  The same geometric
definition given for ordinary polygons works here (generically)
and commutes with projective transformations.

In the next section we will describe coordinates on $\cP_n$.
These coordinates identify $\cP_n$ as an open dense
subset of $\R^{2n}$.  Sometimes we will simply
identify $\cP^n$ with $\R^{2n}$.  
The space $\cC_n$ is much more complicated;  it is an
open dense subset of a codimension $8$ subvariety of
$\R^{2n}$.

\begin{remark}
{\rm
If $n\not=3m$, then it seems useful to impose the simple
condition that $v_i,v_{i+1},v_{i+2}$ are in general
position for all $i$.  With this condition,
$\cP_n$ is isomorphic to the space of 
difference equations of the form
\begin{equation} 
\label{recur}
V_{i}=a_i \,V_{i-1}-b_i\,V_{i-2} + V_{i-3},
\end{equation}
where $a_i,b_i\in\C$ or $\R$ are $n$-periodic:
$a_{i+n}=a_i$ and $b_{i+n}=b_i$, for all $i$.
Therefore, $\cP_n$ is just a $2n$-dimensional vector space, provided $n\not=3m$.
Let us also mention that the spectral theory of difference operators of type
(\ref{recur}) is a classical domain (see \cite{Kri} and references therein).
}
\end{remark}

\subsection{The corner coordinates}\label{coord}
\label{CORNERDEF}
\label{LEFT}

Following \cite{Sch3}, we define local coordinates $(x_1,\ldots,x_{2n})$ on the space $\cP_n$
and give the explicit formula for the pentagram map.  

Recall that the (inverse) \textit{cross ratio} of $4$ collinear points
in $\pP^2$ is given by
\begin{equation}
\label{ICR}
[t_1,t_2,t_3,t_4]=
\frac{(t_1-t_2)\,(t_3-t_4)}{(t_1-t_3)\,(t_2-t_4)},
\end{equation}
where $t$ is (an arbitrary) affine parameter.

\begin{figure}[hbtp]
\centering
\includegraphics[height=1.8in]{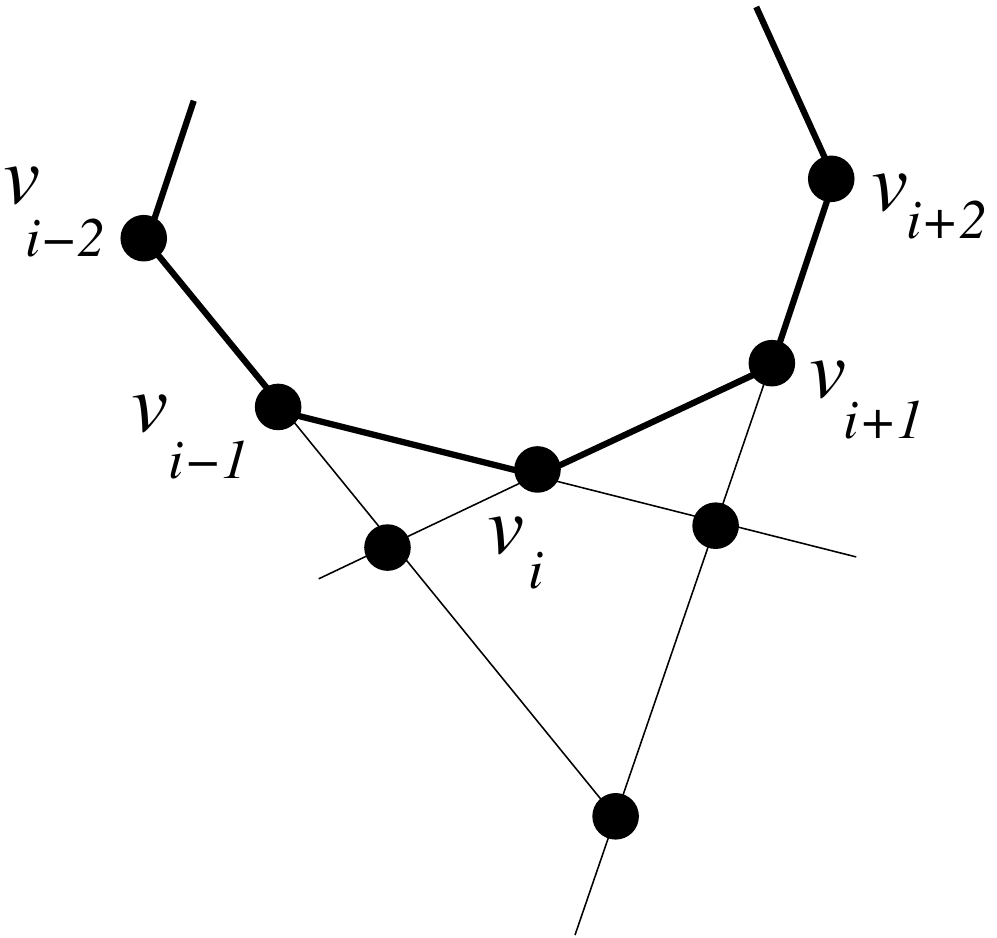}
\newline
\caption{Definition of the corned invariants}
\label{Fig2}
\end{figure}

We define
\begin{equation}
\label{CoD}
\begin{array}{rcl}
x_{2i-1}&=&
\displaystyle
\left[
v_{i-2},\,v_{i-1},\,
\left(
(v_{i-2},v_{i-1})\cap(v_i,v_{i+1})
\right),\,
\left(
(v_{i-2},v_{i-1})\cap(v_{i+1},v_{i+2})
\right)
\right]\\[10pt]
x_{2i+0}&=&
\displaystyle
\left[
v_{i+2},\,v_{i+1},\,
\left(
(v_{i+2},v_{i+1})\cap(v_i,v_{i-1})
\right),\,
\left(
(v_{i+2},v_{i+1})\cap(v_{i-1},v_{i-2})
\right)
\right]\\[10pt]

\end{array}
\end{equation}
where $(v,w)$ stands for the line through $v,w\in\pP^2$,
see Figure \ref{Fig2}.
The functions $(x_1,\ldots,x_{2n})$ are cyclically ordered:
$x_{i+2n}=x_i$.  They provide a system of local coordinates
on the space $\cP_n$ called the \textit{corner invariants}, cf.~\cite{Sch3}.

\begin{remark}
{\rm
a) The index $2i+0$ just means $2i$.  The zero is present to
align the equations. 

b) The right hand side of the second equation is obtained from
the right hand side of the first equation just by swapping the
roles played by $(+)$ and $(-)$.  In light of this fact, it
might seem more natural to label the variables so that the
second equation defines $x_{2i+1}$ rather than $x_{2i+0}$.
The corner invariants would then be indexed by odd integers.
In Section \ref{TAN} we will present an alternate labelling
scheme which makes the indices work out better.

c) Continuing in the same vein, we remark that
there are two useful ways to label the corner
invariants.  In \cite{Sch3} one uses the variables
$x_1,x_2,x_3,x_4,...$ whereas in  \cite{OST,ST2} one uses
the variables $x_1,y_1,x_2,y_2,...$.  The explicit
correspondence between the two labeling schemes is
$x_{2i-1}\to{}x_i,\,x_{2i}\to{}y_i$.  We call the former
convention the {\it flag convention\/} whereas we call the
latter convention the {\it vertex convention\/}.  The reason
for the names is that the variables $x_1,x_2,x_3,x_4$
naturally correspond to the flags of a polygon, as
we will see in Section \ref{TAN}.  The variables $x_i,y_i$
correspond to the two flags incident to the $i$th vertex.
}
\end{remark}

Let us give an explicit formula for the pentagram map in the corner coordinates.
Following~\cite{OST}, we will choose the \textit{right} labelling\footnote{
To avoid this choice between the left or right labelling one can consider the square
$T^2$ of the pentagram map.
.} 
of the vertices
of $T(P)$, see Figure \ref{choices}.
One then has (see  \cite{Sch3}):

\begin{figure}[hbtp]
\centering
\includegraphics[height=130pt]{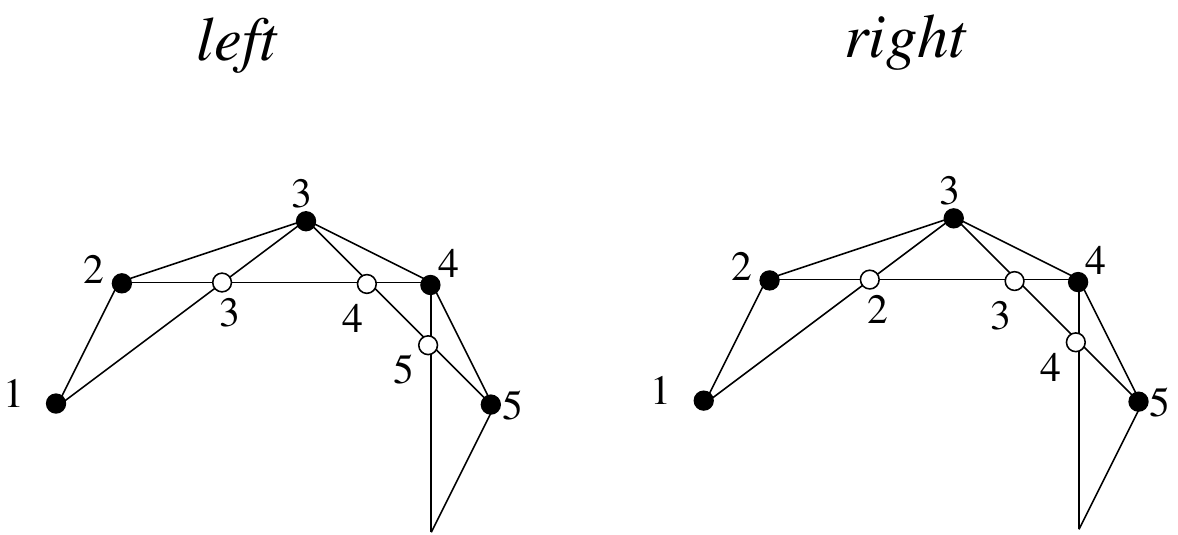}
\newline
\caption{Left and right labelling}
\label{choices}
\end{figure}
\begin{equation}
\label{ExpXEq}
T^*x_{2i-1}=x_{2i-1}\,\frac{1-x_{2i-3}\,x_{2i-2}}{1-x_{2i+1}\,x_{2i+2}},
\qquad 
T^*x_{2i}=x_{2i+2}\,\frac{1-x_{2i+3}\,x_{2i+4}}{1-x_{2i-1}\,x_{2i}},
\end{equation}
where $T^*x_i$ stands for the pull-back of the coordinate functions.

\subsection{Rescaling and the spectral parameter}\label{SpeC}

Equation (\ref{ExpXEq}) has an immediate consequence:
a \textit{scaling symmetry} of the pentagram map.

Consider a one-parameter group $\R^*$ (or~$\C^*$ in the complex case)
acting on the space $\cP_n$ multiplying the coordinates
by $s$ or $s^{-1}$  according to parity:
\begin{equation}
\label{rescal}
\textstyle
R_t:\,(x_1,x_2,x_3\ldots,x_{2n})\to
(s\,x_1,\,s^{-1}\,x_2,\,s\,x_3,\,\ldots,\,s^{-1}x_{2n}).
\end{equation}
It follows from (\ref{ExpXEq}),
that the pentagram map commutes with the rescaling operation.

We will call the parameter $s$ of the rescaling symmetry the \textit{spectral parameter} since it 
defines a one-parameter deformation of the monodromy, $M_s$.
Note that the notion of spectral parameter is extremely useful in the theory
of integrable systems.

\subsection{The Poisson bracket}\label{PoS}

Recall that a \textit{Poisson bracket} on a manifold is a Lie bracket $\{.,.\}$
on the space of functions satisfying the Leibniz rule:
$$
\{F,GH\}=\{F,G\}H+G\{F,H\},
$$
for all functions $F,G$ and $H$.
The Poisson bracket is an essential ingredient of the Liouville-Arnold integrability \cite{Arn}.

Define the following Poisson structure
on~$\cP_n$.
For the coordinate functions we set
\begin{equation}
\label{PoBr}
\{x_i,\,x_{i+2}\}=(-1)^i\,x_i\,x_{i+2},
\end{equation}
and all other brackets vanish.
In other words, the Poisson bracket $\{x_i,x_j\}$ of two coordinate functions
is different from zero if and only if $|i-j|=2$.
The Leibniz rule then allows one to extend  the Poisson bracket to all polynomial
(and rational) functions.

Note that the Jacobi identity obviously  holds.
Indeed, the bracket (\ref{PoBr}) has constant coefficients
when considered in the logarithmic coordinates $\log{}x_i$.

\begin{proposition}
\label{PreProp}
The pentagram map preserves the Poisson bracket (\ref{PoBr}).
\end{proposition}

\proof
This is an easy consequence of formula (\ref{ExpXEq}),
see \cite{OST} (Lemma 2.9), for the details.
\proofend

Recall that a Poisson structure is a way to associate a vector field to a function.
Given a function $f$ on $\cP_n$, the corresponding vector field $X_f$ is called
the \textit{Hamiltonian vector field} defined by $X_f(g)=\{f,g\}$
for every function $g$.
In the case of the bracket (\ref{PoBr}), the explicit formula is as follows:
\begin{equation}
\label{Ham}
X_f=\sum_{i-j=2}
(-1)^{\frac{i+j}{2}}\,
\,x_i\,x_j\left(
\frac{\partial{}f}{\partial{}x_i}\,\frac{\partial}{\partial{}x_j}
-\frac{\partial{}f}{\partial{}x_j}\,\frac{\partial}{\partial{}x_i}
\right).
\end{equation}
Note that the definitions of the Poisson structure in terms of
the bracket of coordinate functions (\ref{PoBr}) and in terms of the Hamiltonian vector fields (\ref{Ham}) are equivalent.

Geometrically speaking, Hamiltonian vector fields are  defined as the image of the map
\begin{equation}
\label{HaMap}
X: T^*_x\,\cP_n\to{}T_x\,\cP_n
\end{equation}
at arbitrary point $x\in\cP_n$.
The kernel of $X$ at a generic point is spanned by the differentials of the Casimir functions, that is, the functions that Poisson commute with all functions.

\begin{remark}
{\rm
The cluster algebra approach of \cite{Gli} also provides a Poisson bracket, 
invariant with respect to the pentagram map (see the book \cite{GSV}).
It can be checked that this cluster Poisson bracket is induced  by the bracket (\ref{PoBr}).
}
\end{remark}

\subsection{The rank of the Poisson bracket and the Casimir functions}\label{PoC1}

The corank of a Poisson structure is the dimension of the kernel of the map $X$ in (\ref{HaMap}), that is, the dimension of the space generated by the differentials of the Casimir functions.

\begin{proposition}
\label{CasProp}
The Poisson bracket (\ref{PoBr}) has corank $2$ if $n$ is odd and corank $4$ is $n$ is even;
the functions
\begin{equation}
\label{casimir0}
O_n=x_1x_3\cdots{}x_{2n-1},
\qquad
E_n=x_2x_4\cdots{}x_{2n}
\end{equation}
for arbitrary $n$ and the functions
\begin{equation}
\label{casimir1}
O_{\frac{n}{2}}=
\prod_{1\leq{}i \leq\frac{n}{2}}x_{4i-1}+
\prod_{1\leq{}i \leq\frac{n}{2}} x_{4i+1},
\qquad
E_{\frac{n}{2}}=
\prod_{1\leq{}i \leq\frac{n}{2}}x_{4i} +
\prod_{1\leq{}i \leq\frac{n}{2}} x_{4i+2},
\end{equation}
for even $n$,
are the Casimirs of the Poisson bracket (\ref{PoBr}).
\end{proposition}

\proof
First, one checks that the functions (\ref{casimir0}) and (\ref{casimir1}) are indeed Casimir functions
(for arbitrary $n$ and for even $n$, respectively).
To this end, it suffices to consider the brackets of (\ref{casimir0}) and (\ref{casimir1}), if $n$ is even, with the coordinate functions $x_i$.

Second, one checks that the corank of the Poisson bracket is
equal to $2$, for odd $n$ and $4$, for even $n$.
The corank is easily calculated in the coordinates $\log{}x_i$,
see \cite{OST}, Section 2.6 for the details.
\proofend

It follows that the Casimir functions are of the form $F \left(O_n,E_n \right)$, if $n$ is odd,
and of the form $F(O_{n/2},E_{n/2},O_n,E_n)$, if $n$ is even.  In both cases
the generic symplectic leaves of the Poisson structure
have dimension $4[(n-1)/2]$.

\begin{remark}
{\rm
If $n$ is even, then the Casimir functions can be written in a more simple manner:
$$
\Big\{
\prod_{1\leq{}i \leq\frac{n}{2}}x_{4i-1},
\quad
\prod_{1\leq{}i \leq\frac{n}{2}} x_{4i+1},
\quad
\prod_{1\leq{}i \leq\frac{n}{2}}x_{4i},
\quad
\prod_{1\leq{}i \leq\frac{n}{2}} x_{4i+2}
\Big\}
$$
instead of (\ref{casimir0}) and (\ref{casimir1}).
}
\end{remark}

\subsection{Two constructions of the monodromy invariants}\label{MoIn}

The second main ingredient of the Liouville-Arnold theory is a set of
Poisson-commuting invariant functions.
In this section, we recall the construction \cite{Sch3} of a set of first integrals
of the pentagram map
$$
O_1,\ldots,O_{\left[\frac{n}{2}\right]},O_n,\;
E_1,\ldots,E_{\left[\frac{n}{2}\right]},E_n
$$
called the \textit{monodromy invariants}.
In other words, we will define $n+1$ invariant function on~$\cP_n$, if $n$ is odd,
and $n+2$ invariant function on $\cP_n$, if $n$ is even.
The monodromy invariants are polynomial in the coordinates (\ref{CoD}).
Algebraic independence of these polynomials was proved in \cite{Sch3}.
Note that $O_n$ and $E_n$ are the Casimir functions (\ref{casimir0}) and, for even $n$,
the functions $O_{\frac{n}{2}}$ and $E_{\frac{n}{2}}$ are as in (\ref{casimir1}).

The indexing of the function $O_i,E_j$ corresponds to their \textit{weight}.
More precisely, we define the weight of the coordinate functions by
\begin{equation}
\label{Weight}
\left|x_{2i+1}\right|=1,
\qquad
\left|x_{2i}\right|=-1.
\end{equation}
Then, $|O_k|=k$ and $|E_k|=-k$.
We give two definitions of the monodromy invariants.
In \cite{Sch3} it is proved that the two definitions
are equivalent.

\bigskip

{\bf A}.
{\it The geometric definition}.
Given a twisted $n$-gon (\ref{PoT}), the corresponding monodromy has a unique lift
to $\SL_3$. By slightly abusing notation, we again denote this
matrix by $M$.  
The two traces, $\mathrm{tr}(M)$ and $\mathrm{tr}(M^{-1})$,
are preserved by the pentagram map (this is a consequence of the projective invariance of $T$).
These traces are rational functions in the corner invariants.
Consider the following two functions:
$$
\widetilde\Omega_1=\mathrm{tr}(M)\,O_n^{\frac{2}{3}}\,E_n^{\frac{1}{3}},
\qquad
\widetilde\Omega_2=\mathrm{tr}(M^{-1})\,O_n^{\frac{1}{3}}\,E_n^{\frac{2}{3}}.
$$
It turns out that
$\widetilde \Omega_1$ and $\widetilde \Omega_2$ are 
polynomials in the corner invariants (see \cite{Sch3}).
Since the pentagram map preserves the monodromy, and 
$O_n$ and $E_n$ are invariants, the two functions
$\widetilde \Omega_1$ and $\widetilde \Omega_2$
are also invariants.
We then have:
\begin{equation}
\label{OandE}
\widetilde \Omega_1=\sum_{k=0}^{[n/2]} O_k,
\qquad
\widetilde \Omega_2=\sum_{k=0}^{[n/2]} E_k,
\end{equation}
where $O_k$ has weight $k$ and $E_k$ has weight $-k$
and where we set
$$
O_0=E_0=1,
$$
for the sake of convenience.
The pentagram map preserves each homogeneous component individually
because it commutes with the rescaling (\ref{rescal}).

Notice also that, if $n$ is even, then $O_\frac{n}{2}$ and $E_\frac{n}{2}$
are precisely the Casimir functions (\ref{casimir1}).
However, the invariants $O_n$ and $E_n$ do not enter the formula (\ref{OandE}).

\bigskip

{\bf B}.
{\it The combinatorial definition}.
Together with the coordinate functions $x_i$, we consider the following ``elementary monomials''
\begin{equation}
\label{uniT}
X_i:=x_{i-1}\,x_{i}\,x_{i+1},
\qquad i=1,\ldots,2n.
\end{equation}

Let $O(X,x)$ be a monomial of the form
$$
O=X_{i_1}\cdots{}X_{i_s}\,x_{j_1}\cdots{}x_{j_t},
$$
where $i_1,\ldots,i_s$ are even and $j_1,\ldots,j_t$ are odd.
Such a monomial is called \textit{admissible} if the Poisson brackets
$\{X_{i_r},X_{i_u}\}$ and $\{X_{i_r},x_{j_u}\}$ and $\{x_{j_r},x_{j_u}\}$
of all the elementary monomials entering~$O$ vanish.

The weight of the above monomial is
$$
|O|=s+t,
$$
see (\ref{Weight}).
For every admissible monomial, we also define the \textit{sign} of $O$ via
$$
\mathrm{sign}(O):=(-1)^t.
$$
The invariant $O_k$ is defined as the alternated sum of all the admissible
monomials of weight $k$:
\begin{equation}
\label{InvOk}
O_k=\sum_{|O|=k}\mathrm{sign}(O)\,O,
\qquad
k\in\left\{1,2,\ldots,\left[\frac{n}{2}\right]\right\}.
\end{equation}
It is proved in \cite{Sch3} that this definition of $O_k$ coincides with (\ref{OandE}).

\begin{example}
{\rm
The first two invariants are:
$$
O_1=\sum_{i=1}^n
\left(
X_{2i}-x_{2i+1}
\right),
\qquad
O_2=\sum_{|i-j|\geq2}
\left(
x_{2i+1}\,x_{2j+1}-
X_{2i}\,x_{2j+1}+X_{2i}\,X_{2j+2}
\right),
$$
for $n\leq5$ the above formulas simplify, see \cite{OST}.
}
\end{example}

The definition of the functions $E_k$ is exactly the same, except that
the roles of {\it even\/} and {\it odd\/} are swapped.

\begin{remark} 
{\rm There is an elegant way to
define the monodromy invariants in terms of determinants.
See \cite{ST2}.}
\end{remark}

\subsection{The monodromy invariants Poisson commute}\label{PoC2}

In this section we give a complete proof of the following
result, which was claimed in \cite{OST}.

\begin{theorem}
\label{PoCThm}
The monodromy invariants Poisson commute with each other, i.e.,
$$
\{O_i,O_j\}=\{O_i,E_j\}=\{E_i,E_j\}=0,
$$
for all $i,j$ indexing the monodromy invariants.
Hence, the Hamiltonian vector fields corresponding to the monodromy
invariants $X_{O_i},\,X_{E_i}$ commute with each other.
\end{theorem}

\proof
The second statement is a consequence of the
first statement.  So, we will just prove the first statement.

We begin with a prelimianry discussion of
how the Poisson bracket interacts with the
elementary monomials defined above.
The Poisson brackets of elementary monomials
\begin{equation}
\label{Element1}
\left\{X_i,X_{i+2}\right\}=(-1)^{i+1}\,X_iX_{i+2},
\qquad
\left\{X_i,X_{i+4} \right\}=(-1)^{i+1}\,X_iX_{i+4},
\end{equation}
together with
\begin{equation}
\label{Element2}
\left\{x_i,X_j \right\}=
\left\{
\begin{array}{rl}
(-1)^i\,x_iX_j, & j=i+1,\,i+2,\,i+3,\\[8pt]
(-1)^{i+1}\,x_iX_j, & j=i-3,\,i-2,\,i-1,
\end{array}
\right.
\end{equation}
immediately follow from the definition (\ref{PoBr}).
All  other brackets $\{X_i,X_j\}$, as well as~$\{X_i,x_j\}$, vanish.
\newline

Now we are ready for the main argument.
Consider first the Poisson bracket $\{O_k,O_m\}$.
This is a sum of the monomials of the form
$$
m=X_{i_1}\cdots{}X_{i_s}x_{j_1}\cdots{}x_{j_t},
$$
where $i_1,\ldots,i_s$ are even and $j_1,\ldots,j_t$ are odd.
Indeed, by definition of the Poisson structure~(\ref{PoBr}), the bracket
of two monomials is proportional to their product, so that
the above bracket contains only the monomials entering $O_k$ and $O_m$.

The monomial $m$ is not necessarily admissible.
There can be squares (some $i$'s or $j$'s may coincide), but no cubes
or higher degrees.
We want to prove that the numeric coefficient of every such monomial in $\{O_k,O_m\}$ is zero.

We define an \textit{oriented graph} with the set of vertices
$\{X_{i_1},\ldots,X_{i_s},\,x_{j_1},\ldots,x_{j_t}\}$ corresponding to the elementary monomials
in $m$; the oriented arrows joining the
vertices whenever their Poisson bracket is different from zero, the orientation being given
by the sign of the bracket.
Recall that all the non-zero brackets of elementary monomials
are listed in (\ref{Element1}) and (\ref{Element2}).

\begin{lemma}
\label{Single}
If two indices coincide, $i_r=i_u$ or $j_r=j_u$, then the
corresponding connected component of the graph consists of one element.
\end{lemma}

\proof
If $i_r=i_u$, then $X_{i_r}=X_{i_u}$ belongs both to $O_i$ and $O_j$.
By the admissibility condition, this implies all the Poisson brackets of $X_{i_r}$
with the other elementary monomials from $m$ vanish.
\proofend

The above lemma allows one to assume that all the indices in $m$ are different:
$i_r\not=i_u$ and $j_r\not=j_u$.

\begin{lemma}
\label{GraphL}
The above defined graph has

(i)
no $3$-cycles;

(ii)
no vertices with more than one outgoing or ingoing arrows;
in other words, the graph does not have the following vertices:
$$
\begin{CD}
\;@<<< a @> >> \;
\qquad
\;@>>> a @<<< \;
\end{CD}
$$
where $a=X_{i_r}$ or $x_{j_u}$.
\end{lemma}

\proof
(i) Assume there is a $3$-cycle.
Then at least two of the corresponding elementary monomials
belong  to the decomposition of either $O_i$ or  $O_j$.
The monomials are joined by an arrow, thus their Poisson bracket does not vanish.
This leads to a contradiction since all the
monomials in $O_i$ are admissible (see Section \ref{MoIn}, definition B).

(ii) To show that no vertex of the graph
can have more than one outgoing or ingoing arrows, one has to analyze
formulas  (\ref{Element1}) and (\ref{Element2}).
Since $i_r$ are even and $j_u$ are odd, a vertex~$X_{i_r}$ can be joined by an outgoing arrow
to the following vertices (provided they belong to the graph):
$X_{i_r-2},\,X_{i_r-4},\,x_{i_r+1},x_{i_r+3}$.
In all of these cases, we obtain a 3-cycle, which is a contradiction to part (i) of the lemma.
\proofend

The above lemma implies the following statement.
\begin{corollary}
\label{GraphC}
The graph has no branching (i.e., vertices with three or more adjacent arrows).
\end{corollary}
\noindent
Indeed, a branching point has more than one out- or ingoing arrows:
$$
\begin{CD}
@.@AAA  \\
@> >>a@> >>
\end{CD}
\qquad
\begin{CD}
@.@AAA  \\
@> >>a@<<<
\end{CD}
$$

\begin{remark}
{\rm
One can also show that the constructed graph has no $k$-cycles for arbitrary $k$, that is,
every connected component of the graph 
is of type $A_k$ oriented in the standard way:
$$
\begin{CD}
a_1@> >>a_2@> >> \cdots @> >> a_k
\end{CD}
$$
but we will not use this in the proof.
}
\end{remark}

Let us finally deduce $\{O_i,O_j\}=0$ from Lemma \ref{GraphL} and Corollary \ref{GraphC}. 
Every element $X_{i_r}$ and~$x_{i_u}$ in the monomial $m$ belongs either to $O_i$, or to $O_j$. 
If the constructed graph contains at least three elements, then is has a fragment:
$$
\begin{CD}
a_1@> >>a_2@> >> a_3
\end{CD}
$$
where either $a_1,a_3\in O_i$ and $a_2\in O_j$ or the other way around.
It follows from the Leibniz identity that the element $a_2$ contributes twice
 in $\{O_i,O_j\}$, namely in
$\{a_1,a_2\}$ and in $\{a_3,a_2\}$, with the opposite signs.

We proved $\{O_i,O_j\}=0$, except the case where the graph is
of type $A_2$, i.e., contains only two elements:
$$
\begin{CD}
a_1@> >>a_2
\end{CD}
$$
with, say, $a_1\in O_i, a_2\in O_j$.
But in this last case,  the Poisson brackets of the elementary monomials $a_1$ and $a_2$
with all the other elementary monomials in $m$ vanish.
By construction of the invariants, $O_i$ and $O_j$ are symmetric
with respect to the monomials $a_1$ and $a_2$.
It follows that the monomial $m$ appears twice in $\{O_i,O_j\}$,
with the opposite signs.
This completes the proof that $\{O_i,O_j\}=0$.

The proof of $\{E_i,E_j\}=0$ is identically the same (with odd and even indices exchanged).
It remains to consider the bracket $\{O_i,E_j\}$.

We will apply the same idea and construct a graph for every monomial in $\{O_i,E_j\}$.
Recall that $E_j$ contains the admissible monomials 
$E=X_{i_1}\cdots{}X_{i_s}\,x_{j_1}\cdots{}x_{j_t},$
where the indices $i_1,\ldots,i_s$ are odd and $j_1,\ldots,j_t$ are even.
Analyzing the brackets  (\ref{Element1}) and (\ref{Element2}),
we see that the graph corresponding to any monomial in $\{O_i,E_j\}$
is of the form
$$
\begin{CD}
\cdots@> >>x_i@> >>X_{i+2}@> >> x_{i+4}@> >>X_{i+6}@> >>\cdots
\end{CD}
$$
and the $X$'s and $x$'s belong to the different functions.
We observe that
$$
\{x_i,\,X_{i+2}\}=-\{X_{i+1},\,x_{i+3}\},
$$
and if $x_i\in{}O_i$ and $X_{i+2}\in{}E_j$ then
$O_i$ and $E_j$ are symmetric with respect to the exchange of
$x_i$ with $X_{i+1}$ and of $X_{i+2}$ with $x_{i+3}$,
respectively.
The monomial $m$ appears twice with the opposite signs.
This completes the proof of Theorem \ref{PoCThm}.
\proofend

In \cite{Sch3} it is proved that the monodromy invariants are algebraically
independent.  The argument is rather complicated, but it is very similar
in spirit to the related independence proof we give in Section \ref{INDEP}.
The algebraic independence result combines with
Theorem \ref{PoCThm} to establish the
integrability of the pentagram map on the space $\cP_n$.
Indeed, the Poisson bracket (\ref{PoBr}) defines a symplectic foliation on $\cP_n$, the
symplectic leaves being locally described as levels of the Casimir functions, see Proposition \ref{CasProp}.
The number of the remaining invariants is exactly half of the dimension of the symplectic leaves.
The classical Liouville-Arnold theorem~\cite{Arn} is then applied.

\section{Integrability on $\cC_n$ modulo a calculation}

The general plan of the proof of Theorem \ref{Main} is as follows.

\begin{enumerate}
\item
We show that the Hamiltonian vector fields on $\cP_n$ corresponding
to the monodromy invariants are \textit{tangent} to the subspace $\cC_n$,

\item
We restrict the monodromy invariants to $\cC_n$
and show that the dimension of a generic level set is $n-4$ if $n$ is odd
and $n-5$ if $n$ is even.

\item
We show that there are exactly the same number
of independent Hamiltonian vector fields.
\end{enumerate}

In this section, we prove the first statement and also show that
the dimension of the level sets is \textit{at most} $n-4$ if $n$ is odd
and $n-5$ if $n$ is even, and similarly for the number
of independent Hamiltonian vector fields.
The final step of the proof that this upper bound is 
actually the lower one will be done in the next two sections. 
This final step is a nontrivial calculation that
comprises the bulk of the paper.

\subsection{The Hamiltonian vector fields are tangent to $\cC_n$}\label{TanS}

The space $\cC_n$ is a subvariety of $\cP_n$ having codimension $8$.
It turns out that one can give explicit equations for this
variety.  See Lemma \ref{variety}.  (These equations do not
play a role in our proof, but they are useful to have.)

The following statement is  essentially a consequence of Theorem \ref{PoCThm}.
This is an important step of the proof of Theorem \ref{Main}.

\begin{proposition}
\label{TanProp}
The Hamiltonian vector field on $\cP_n$ corresponding to a monodromy invariant
is tangent to $\cC_n$.
\end{proposition}

\proof
The space $\cP_n$ is foliated by isomonodromic submanifolds 
that are generically of codimension $2$ and are defined by the condition that the monodromy has fixed eigenvalues. Hence the isomonodromic submanifolds can be defined as the level surfaces of two functions, 
$\mathrm{tr}(M)$ and $\mathrm{tr}(M^{-1})$. 
This foliation is singular, and $\cC_n$ is a singular leaf of codimension $8$. 
We note that the versal deformation of $\cC_n$
is locally isomorphic to $\SL(3)$ partitioned into the conjugacy equivalence classes.

Consider a monodromy invariant, $F\,(=O_i$ or $E_i$), and its Hamiltonian vector field, $X_F$.
We know that the Poisson bracket $\{F,\,\mathrm{tr}(M)\}=0$,
since all monodromy invariants Poisson commute
and $\mathrm{tr}(M)$ is a sum of monodromy invariants. 
Hence $X_F$ is tangent to the generic leaves of the isomonodromic foliation on $\cP_n$. 
Let us show that $X_F$ is tangent to $\cC_n$ as well.

In a nutshell, this follows from the observation that the tangent space to $\cC_n$ at a smooth point $x_0$ is the intersection of the limiting positions of the tangent spaces to the isomonodromic leaves at points $x$ as $x$ tends to $x_0$. Assume then that  $X_F$ is transverse to $\cC_n$ at point $x_0\in \cC_n$. Then $X_F$ will be also transverse to an isomonodromic leaf at some point $x$ close to $x_0$, yielding a contradiction.

More precisely, we can apply a projective transformation so that the vertices $V_1, V_2,V_3,V_4$ of a twisted $n$-gon $V_1,V_2,\dots$ become the vertices of a standard square. This gives a local identification of $\cP_n$ with  the set of tuples $(V_5,\dots,V_n; M)$ where $M$ is the monodromy, the projective transformation that takes the quadruple 
$(V_1,V_2,V_3,V_4)$ to $(V_{n+1}, V_{n+2}, V_{n+3}, V_{n+4}).$
The space of closed $n$-gons is characterized by the condition that $M$ is the identity. Thus we have locally identified $\cP_n$ with $\cC_n \times \SL(3)$. In particular, we have a projection $\cP_n \to \SL(3)$, and the preimage of the identity is $\cC_n$. The  isomonodromic leaves project to the conjugacy equivalent classes in $\SL(3)$.

Thus our proof reduces to the following fact about the group $\SL(3)$ (which holds for $\SL(n)$ as well). 

\begin{lemma}
\label{limpos}
Consider the singular foliation of $\SL(3)$ by the conjugacy equivalence classes, and let $T_X$ be the tangent space to this foliation at $X\in \SL(3)$. Then the intersection, over all $X$, of the limiting positions of the spaces $T_X$, as $X \to 1\!\!1$, is trivial (here $1\!\!1\in \SL(3)$ is the identity). 
\end{lemma}

\proof Let $B\in\SL(3)$, and let $B+\varepsilon C$ be an
infinitesimal deformation within the conjugacy equivalence class. Then
$$
\mathrm{tr}\left(B+\varepsilon C\right)=
\mathrm{tr}(B), \qquad \mathrm{tr} \left((B+\varepsilon C)^2 \right)=
\mathrm{tr} \left(B^2 \right),
$$
hence $\mathrm{tr}(C)=0$ and $\mathrm{tr}(BC)=0$,
and also $\mathrm{tr}(B^{-1} C)=0$ since $\det (B+\varepsilon C)=1$. 
Thus the tangent space to a conjugacy
equivalent class of $B$ is given by
$$
\mathrm{tr}(C)=\mathrm{tr}(BC)=\mathrm{tr}(B^{-1} C)=0.
$$
Now let $B=1\!\!1+\varepsilon A$, a point in an infinitesimal neighborhood of the identity $1\!\!1$; we have $\mathrm{tr}(A)=0$. 
Then our conditions on~$C$ implies $\mathrm{tr}(C)=\mathrm{tr}(AC)=0$. 
Since $\mathrm{tr}(AC)$ is a non-degenerate quadratic form, 
an element $C\in \mathrm{sl}(3)$ satisfying $\mathrm{tr}(AC)=0$ for all $A\in \mathrm{sl}(3)$
has to be zero.
\proofend

In view of what we said above, this implies the proposition.
\proofend

\subsection{Identities between the monodromy invariants}\label{IdenS}

In this section, we consider the restriction of the monodromy invariants from the space
of all twisted $n$-gons to the space $\cC_n$ of closed $n$-gons.
We show that these restrictions satisfy 5 non-trivial relations,
 whereas their differentials, considered as covectors in $\cP_n$ whose foot-points belong to $\cC_n$, 
satisfy 3 non-trivial relations.  
These relations are also mentioned  in \cite{OST} and \cite{Sol}.
In Sections \ref{ThEnd} and \ref{TAN}, we will prove that there are no other relations
between the monodromy invariants on $\cC_n$ and their differentials along $\cC_n$.

We remark that, strictly speaking, the identities established in this section
are not needed for the proof of our main result.  For the main result, all we need to know is that
there are enough commuting flows to fill out what could be ({\it a priori\/}, with out the
results in this section) a union of level sets of the monodromy invariants.  Thus,
the reader interested only in the main result can skip this section.

\begin{theorem} 
\label{5and3xy}
(i) The restrictions of the monodromy integrals to $\cC_n$ satisfy the following five identities:
\begin{equation}
\label{Rel1}
\begin{array}{rcl}
\displaystyle
\sum_{j=0}^{[n/2]}  O_j = 3\,E_n^{\frac{1}{3}}\,O_n^{\frac{2}{3}},
&&\displaystyle
\sum_{j=0}^{[n/2]} E_j = 3\,E_n^{\frac{2}{3}}\,O_n^{\frac{1}{3}},\\[16pt]
\displaystyle
\sum_{j=1}^{[n/2]} j\, O_j = n\, E_n^{\frac{1}{3}}\,O_n^{\frac{2}{3}},
&&\displaystyle
 \sum_{j=1}^{[n/2]} j\, E_j = n\, E_n^{\frac{2}{3}}\,O_n^{\frac{1}{3}},\\[16pt]
\displaystyle
E_n^{\frac{1}{3}}\, \sum_{j=1}^{[n/2]}  j^2 \,O_j&=&
\displaystyle
 O_n^{\frac{1}{3}}\, \sum_{j=1}^{[n/2]} j^2\, E_j.
\end{array}
\end{equation}
(ii) The differentials of the monodromy integrals along $\cC_n$ satisfy the three identities:
\begin{equation}
\label{Rel2}
\begin{array}{rcl}
\displaystyle
\sum_{j=1}^{[n/2]}  dO_j&=&
\displaystyle
2\, E_n^{\frac{1}{3}}\,O_n^{-\frac{1}{3}}\, dO_n + 
E_n^{-\frac{2}{3}}\,O_n^{\frac{2}{3}} dE_n,\\[16pt]
\displaystyle
\sum_{j=1}^{[n/2]} dE_j&=&
\displaystyle
2\, E_n^{-\frac{1}{3}}\,O_n^{\frac{1}{3}} dE_n + 
E_n^{\frac{2}{3}}\,O_n^{-\frac{2}{3}}\, dO_n, \\[16pt]
\displaystyle
O_n^{\frac{1}{3}}\, \Big(
\sum_{j=1}^{[n/2]} j\ dE_j
\Big) + 
E_n^{\frac{1}{3}} \,\Big(
\sum_{j=1}^{[n/2]}  j\ dO_j
\Big)&=&
\displaystyle
n\, E_n^{\frac{2}{3}}O_n^{\frac{2}{3}} \left(
E_n^{-1} dE_n + O_n^{-1} dO_n
\right).
\end{array}
\end{equation}
\end{theorem}

\proof
Recall that the monodromy invariants $O_j$ are the homogeneous components of the polynomial 
$O_n^{2/3}E_n^{1/3}\,\mathrm{tr}(M)$ with respect to the rescaling (\ref{rescal}),
where $s=e^t$ for convenience. 
Likewise, the monodromy invariants~$E_j$ are homogeneous components of  
$O_n^{1/3}E_n^{2/3}\,\mathrm{tr}(M^{-1})$.
 Recall also that $O_0=E_0=1$.
 
Denote for simplicity $O_n^{1/3}E_n^{2/3}=U,\, O_n^{2/3}E_n^{1/3}=V$. 
Notice that the monodromy matrix $M$ has the unit determinant. Let 
$e^{\lambda_1},\  e^{\lambda_2},\ e^{\lambda_2}$
be the eigenvalues of $M$. 
One has
\begin{equation} 
\label{zero}
\lambda_1+\lambda_2+\lambda_3\equiv 0.
\end{equation}
We consider a one-parameter family of $n$-gons
depending on the rescaling parameter $t$, such that for $t=0$, the $n$-gon belongs to $\cC_n$.
The monodromy $M=M_t$ also depends on $t$ so that we think of $\lambda_i$ as 
functions of the corner coordinates 
$(x_1,\ldots,x_{2n})$ and of $t$. 
For $t=0$, one has: $\lambda_i=0,\, i=1,2,3$ since $M_0=\mathrm{Id}$.

The eigenvalues of $M^{-1}$ are
$e^{-\lambda_1},\  e^{-\lambda_2},\ e^{-\lambda_2}. $
Since the weights of $O_i$ and $E_j$ are $j$ and $-j$ respectively, the definition of the integrals writes as follows:
$$
e^{\frac{nt}{3}}\, V \left(e^{\lambda_1}+ e^{\lambda_2}+ e^{\lambda_2}\right)
=\sum_{j=0}^{[n/2]} e^{tj}\, O_j,
\qquad
 e^{-\frac{nt}{3}}\, U\left(e^{-\lambda_1}+ e^{-\lambda_2}+ e^{-\lambda_2}\right)
 =\sum_{j=0}^{[n/2]} e^{-tj}\, E_j.
$$
which we rewrite as
\begin{equation}
\label{def}
V \left(e^{\lambda_1}+ e^{\lambda_2}+ e^{\lambda_2}\right)=
\sum_{j=0}^{[n/2]} e^{t(j-\frac{n}{3})} \,O_j,\qquad
U\left(e^{-\lambda_1}+ e^{-\lambda_2}+ e^{-\lambda_2}\right)=
\sum_{j=0}^{[n/2]} e^{-t(j-\frac{n}{3})}\, E_j.
\end{equation}

Setting $t=0$ in these formulas yields the first two identities in (\ref{Rel1}).
Next, differentiate these equations in $t$ :
$$
 V\, \sum_{i=1}^3 \lambda_i' \,e^{\lambda_i}
 =\sum_{j=0}^{[n/2]} \left(j-\frac{n}{3}\right) e^{tj}\, O_j,
$$
where $\lambda_i' =d\lambda_i/dt$,
and similarly for $E_j$.
Set $t=0$, then
the left-hand-side vanishes because $\sum \lambda_i'=0$ due to (\ref{zero}).
Hence
$$
\sum_{j=0}^{[n/2]} j \, O_j = \frac{n}{3}\, \sum_{j=0}^{[n/2]}   O_j= n\,V
$$
due to the first identity in (\ref{Rel1}) and similarly for $E_j$. 
One thus obtains the third and the fourth identity in (\ref{Rel1}). 

To obtain the fifth equation in (\ref{Rel1}), differentiate  the equations (\ref{def}) 
with respect to $t$ twice to get
$$
\begin{array}{rcl}
\displaystyle
V\,\Big(\sum_{i=1}^3 (\lambda_i''+\lambda_i'^2)\, e^{\lambda_i}\Big)&=&
\displaystyle
\sum_{j=0}^{[n/2]} \left(j-\frac{n}{3}\right)^2 e^{tj}\, O_j,\\[16pt]
\displaystyle
U\,\Big(\sum_{i=1}^3 \left(-\lambda_i''+\lambda_i'^2\right)\, e^{\lambda_i}\Big)&=&
\displaystyle
\sum_{j=0}^{[n/2]} \left(j-\frac{n}{3}\right)^2 e^{-tj}\, E_j.
\end{array}
$$
Divide the first equality by $V$, the second by $U$, subtract one from another, and set $t=0$:
$$
2 \sum_{i=1}^3 \lambda_i''=V^{-1} \sum_{j=0}^{[n/2]} 
\left(j-\frac{n}{3}\right)^2  O_j - U^{-1} \sum_{j=0}^{[n/2]} \left(j-\frac{n}{3}\right)^2  E_j.
$$
The left hand side vanishes, due to (\ref{zero}), so
\begin{equation}
 \label{useful}
V^{-1}\, \sum_{j=0}^{[n/2]} \left(j-\frac{n}{3}\right)^2  O_j = 
U^{-1} \,\sum_{j=0}^{[n/2]} \left(j-\frac{n}{3}\right)^2  E_j.
\end{equation}
Therefore
$$
\begin{array}{l}
\displaystyle
V^{-1} \sum_{j=0}^{[n/2]} j^2\, O_j - \frac{2n}{3}\, V^{-1} \,
\sum_{j=0}^{[n/2]} j\, O_j + V^{-1}\, \frac{n^2}{9}\,  \sum_{j=0}^{[n/2]}  O_j
=\\[10pt]
\displaystyle
\hskip 2cm 
U^{-1} \sum_{j=0}^{[n/2]} j^2\, E_j - \frac{2n}{3}\, U^{-1} \,
\sum_{j=0}^{[n/2]} j\, E_j + U^{-1}\, \frac{n^2}{9}\,  \sum_{j=0}^{[n/2]}  E_j.
\end{array}
$$
The second and the third terms on the left and the right hand sides are pairwise equal, 
due to the first four identities in (\ref{Rel1}). This implies the fifth identity (\ref{Rel1}).

To prove (\ref{Rel2}), take differentials of  (\ref{def}):
$$
\label{diff1}
\begin{array}{l}
\displaystyle
V \sum_{i=1}^3 e^{\lambda_i}\, d\lambda_i + 
\Big(\sum_{i=1}^3 e^{\lambda_i}\Big)\, dV =\\ [10pt]
\displaystyle
\Big(\sum_{j=0}^{[n/2]} \left(j-\frac{n}{3}\right) e^{t(j-\frac{n}{3})}\, O_j\Big) dt + 
\sum_{j=0}^{[n/2]} e^{t(j-\frac{n}{3})}\, dO_j,
\end{array}
$$
and
$$
\begin{array}{l}
\displaystyle
-U \sum_{i=1}^3 e^{-\lambda_i}\, d\lambda_i + 
\Big(\sum_{i=1}^3 e^{-\lambda_i}\Big)\, dU = \\[10pt]
\displaystyle
-\Big(\sum_{j=0}^{[n/2]} \Big(j-\frac{n}{3}\Big) e^{-t(j-\frac{n}{3})} E_j\Big)\, dt + 
\sum_{j=0}^{[n/2]} e^{-t(j-\frac{n}{3})}\, dE_j.
\end{array}
$$
Set $t=0$: the first terms on the right hand sides vanish due to (\ref{zero}), 
and the first parentheses on the right hand sides vanish due to (\ref{Rel1}). 
We get
$$
\sum_{j=0}^{[n/2]}  dO_j=3\,dV,
\qquad 
\sum_{j=0}^{[n/2]}  dE_j=3\,dU,
$$
the first two identities in (\ref{Rel2}).

Finally, differentiate the above equations with respect to $t$ and set $t=0$ to obtain:
$$
\begin{array}{rcl}
\displaystyle
V\, \sum_{i=1}^3 \lambda_i'\, d\lambda_i + 
V\, \sum_{i=1}^3 d(\lambda_i') + \Big(\sum_{i=1}^3 \lambda_i'   \Big)\, dV&=&
\displaystyle
\Big(\sum_{j=0}^{[n/2]} \left(j-\frac{n}{3}\right)^2  O_j\Big)\, dt + 
\sum_{j=0}^{[n/2]} \left(j-\frac{n}{3}\right) d O_j,\\[14pt]
\displaystyle
U\, \sum_{i=1}^3 \lambda_i'\, d\lambda_i -
U\, \sum_{i=1}^3\, d(\lambda_i') +
 \Big(\sum_{i=1}^3 \lambda_i'   \Big)\, dU &=&
\displaystyle
\Big(\sum_{j=0}^{[n/2]} \left(j-\frac{n}{3} \right)^2  E_j\Big)\, dt -
 \sum_{j=0}^{[n/2]} \left(j-\frac{n}{3}\right) d E_j.
 \end{array}
$$
Once again, the second and the third sums on the left hand sides vanish, due to (\ref{zero}). Divide the first equation by $V$, the second by $U$, and subtract one from another, using (\ref{useful}):
$$
V^{-1}\, \sum_{j=0}^{[n/2]} \left(j-\frac{n}{3}\right) d O_j + 
U^{-1}\, \sum_{j=0}^{[n/2]} \left(j-\frac{n}{3}\right) d E_j=0.
$$
Hence
$$
V^{-1}\, \sum_{j=0}^{[n/2]} j\, d O_j + U^{-1}\, \sum_{j=0}^{[n/2]} j\, d E_j = 
\frac{n}{3}\, \Big(V^{-1} \sum_{j=0}^{[n/2]}  d O_j + U^{-1} \sum_{j=0}^{[n/2]}  d E_j \Big).
$$
Due to the first two identities in (\ref{Rel2}), the right-hand-side equals 
$n\, (O_n^{-1}\, dO_n + E_n^{-1}\, dE_n)$. 
This yields the third identity in (\ref{Rel2}).
Theorem \ref{5and3xy} is proved.
\proofend

\begin{remark}
\label{rmk}
{\rm 
a)
Let ${\cal E}$ be the Euler vector field that generates the scaling. Then 
$$
{\cal E} (O_j)=j\,O_j,
\qquad 
\ {\cal E} (E_j)=-j\,E_j.
$$
If one evaluates the differentials in the identities  (\ref{Rel2}) on ${\cal E}$, 
one obtains the last three identities in (\ref{Rel1}). 
This is a check that (\ref{Rel1}) and (\ref{Rel2}) are consistent with each other.

b)
Equivalently, (\ref{Rel2}) can be rewritten as
$$
\begin{array}{rcl}
\displaystyle
3\, d O_n &=& 2\, E_n^{-\frac{1}{3}}\,O_n^{\frac{1}{3}} \,
\Big( \sum_{j=1}^{[n/2]}  dO_j \Big)- 
E_n^{-\frac{2}{3}}\,O_n^{\frac{2}{3}}\, \Big( \sum_{j=1}^{[n/2]}  dE_j \Big),\\[16pt]
\displaystyle
3\, d E_n &=& 2\, E_n^{\frac{1}{3}}\,O_n^{-\frac{1}{3}} \,
\Big( \sum_{j=1}^{[n/2]} dE_j \Big)- 
E_n^{\frac{2}{3}}\,O_n^{-\frac{2}{3}}\, \Big( \sum_{j=1}^{[n/2]}  dO_j \Big),\\[12pt]
0&=&
\displaystyle
O_n^{\frac{1}{3}} \,
\Big(3 \sum_{j=1}^{[n/2]} j\ dE_j -n \sum_{j=1}^{[n/2]} dE_j  \Big) +
 E_n^{\frac{1}{3}} \,\Big(3 \sum_{j=1}^{[n/2]}  j\, dO_j -n \sum_{j=1}^{[n/2]}   dO_j  \Big).
 \end{array}
$$

c)
The identities  (\ref{Rel1}) and (\ref{Rel2}) are satisfied in a larger subspace than $\cC_n$,
 consisting of twisted polygons whose monodromy has \textit{equal eigenvalues}.
This subspace has codimension 2 in~$\cP_n$. 

d) In both cases, $n$ odd and $n$ even, the kernel of the Poisson map $X$ (\ref{HaMap}) (spanned by the differentials of the Casimir functions) has zero intersection with the subspace of $T^*\cP_n$ spanned by the relations \ref{Rel2}.
}
\end{remark}

\subsection{Reducing the proof to a one-point computation}\label{ReS}

For ease of exposition, we will give our proof only in
the odd case, and we set $n \geq 7$ odd.  Modulo
changing some of the indices, the even case is similar.
We will explain everything in terms of the odd case and,
at the end of this section, briefly explain what happens
in the even case.

Let ${\cal M\/}$ denote the algebra generated by
the monodromy invariants.   In the Section \ref{ThEnd} we make
the following calculations.

\begin{enumerate}
\item There exist elements $F_1,...,F_{n-2} \in \cal M$
and a point $p \in \cC_n$ such that the differentials
$dF,...,dF_{n-2}$ are linearly independent at $p$.
Therefore,
$dF,...,dF_{n-2}$ are linearly independent at almost all
$q \in \cC_n$. 
\item There exists elements $G_1,...,G_{n-4} \in \cal M$
and a point $p \in \cC_n$ such that the differentials
$dG_1|_{T_p \cC_n},...,dG_{n-4}|_{T_p \cC_n}$ are linearly independent.
Therefore, $dG_1|_{T_q \cC_n},...,dG_{n-4}|_{T_q \cC_n}$ are linearly independent
at almost all $q \in \cC_n$. 
\end{enumerate}
In Calculation 1,
we are computing the differentials
on the ambient space ${\cal P\/}_n$ but evaluating them
at a point of $\cC_n$. In Calculation 2, we are
computing the differentials on the ambient space,
evaluating them at a point of $\cC_n$, {\it and\/}
restricting the resulting linear functionals to the
tangent space of $\cC_n$.   In both calculations,
we are actually evaluating at points in
$\C_n^0$.  In each case, what allows us to make
a conclusion about generic points is that the
monodromy invariants are algebraic.

Calculation 2 combines with Theorem~\ref{5and3xy} to show that
there are exactly $n-4$ algebraically independent monodromy
invariants, when restricted to $\cC_n$.
Hence, the generic common level set of the monodromy
invariants $O_i,E_i$, restricted to $\cC_n$,
has dimension $n-4$.

Next, we wish to prove that these level sets
have locally free action of the abelian group~$\R^{d}$ (or~$\C^{d}$ in the complex case).
For $F \in \cal M$, the Hamiltonian vector field $X_F$ is tangent to $\cC_n$,
by Proposition \ref{TanProp}, and also tangent to the common level set
of functions in $\cal M$.  Finally, by Theorem~\ref{PoCThm},
the Hamiltonian vector fields all commute with each other
(i.e., define an action of the Abelian Lie algebra).
The following lemma finishes our proof.

\begin{lemma}
The Hamiltonian vector fields of the monodromy invariants 
generically span the monodromy level sets on $\cC_n$.
\end{lemma}

\startproof
Let $\wedge^1\cP_n$ denote the space of
$1$-forms on $\cP^n$.  Let $\cal X$ denote the space
of vector fields on $\cC_n$.
Let $d{\cal M\/} \subset  \wedge^1\cP_n$
denote the image of $\cal M$ under the $d$-operator.
Calculation 1 shows that the vector space $d{\cal M\/}$ 
generically has dimension $n-2$ when evaluated at 
points of $\cC^n$.  At the same time, we have the Poisson 
map $X: d{\cal M\/} \to {\cal X\/}$, given by
$$
X(dF)=X_F,
$$
see (\ref{HaMap}).
In the odd
case, the map $X$ has $2$ dimensional kernel, see Remark \ref{rmk} d).  Hence,
$X$ has $n-4$ dimensional image, as desired.
\proofend

Now we explain explicitly how the results above
give us the quasi-periodic motion in the case of
closed convex polygons.  We know from the work in
\cite{Sch1} that the monodromy level sets on
$\cC_n^0$ are compact.  By Sard's Theorem, and by the
calculations above, almost every level set is
a smooth compact manifold of dimension $m=n-4$.  By
Sard's Theorem again, and by the dimension count above,
almost every level set $L$ possesses a framing by 
Hamiltonian vector fields.  That is, there are
$m$ Hamiltonian vector fields on $L$ which are
linearly independent at each point and which define
commuting flows.  These vector fields define
local coordinate charts from $L$ into $\R^m$,
such that the overlap functions are translations.
Therefore $L$ is a finite union of affine
$m$-dimensional tori.  The whole structure
is invariant under the pentagram map, and so
the pentagram map is a translation of $L$
relative to the affine structure on $L$.  This is the
quasi-periodic motion.  Even more explicitly,
some finite power of the pentagram map preserves
each connected component of $L$ and is a constant
shift on each connected component.
\newline
\newline
\noindent
{\bf The Even Case:\/} 
In the even case, we have the following calculations:
\begin{enumerate}
\item There exist elements $F_1,...,F_{n-1} \in \cal M$
and a point $p \in \cC_n$ such that the differentials
$dF,...,dF_{n-1}$ are linearly independent at $p$.
Therefore,
$dF,...,dF_{n-1}$ are linearly independent at almost all
$q \in \cC_n$. 
\item There exists elements $G_1,...,G_{n-3} \in \cal M$
and a point $p \in \cC_n$ such that the differentials
$dG_1|_{T_p \cC_n},...,dG_{n-3}|_{T_p \cC_n}$ are linearly independent.
Therefore, $dG_1|_{T_p \cC_n},...,dG_{n-3}|_{T_p \cC_n}$ are linearly independent
at almost all $q \in \cC_n$. 
\end{enumerate}

In this case, the common level sets generically
have dimension $n-5$ and, again, the Hamiltonian
vector fields generically span these level sets. The situation is summarized in the following table.

\begin{displaymath}
  \begin{array}{c|c|c|c|}
&\hfill{\rm Invariants}& \hfill{\rm Casimirs} & 
\hfill{\rm Level\ sets}\,/\, \hfill{\rm Hamiltonian}\, \hfill{\rm fields}\\
\cline{1-1} \cline{2-2}\cline{3-3}\cline{4-4}
n\, \hfill{\rm odd}&n+1& 2 & d=n-4  \\
\cline{1-1}\cline{2-2} \cline{3-3}\cline{4-4}
n\, \hfill{\rm even}& n+2 & 4 & d=n-5   \\
\cline{1-1}\cline{2-2} \cline{3-3}\cline{4-4}
\end{array}
\end{displaymath}

\section{The linear independence calculation}\label{ThEnd}
\label{INDEP}

\subsection{Overview}

For any given (smallish) value of $n$, one can make the calculations
directly, at a random point, and see that it works.  The
difficulty is that we need to make one calculation for each $n$.
One might say that the idea behind our calculations is
tropicalization.  The monodromy invariants and their
gradients are polynomials with an enormous number
of terms.  We only need to make our calculation at
one point, but we will consider a $1$-parameter family
of points, depending on a parameter $u$.  As $u \to 0$,
the different variables tend to $0$ at different rates. This
sets up a kind of hierarchy (or filtration) on the the monomials
comprising the polynomials of interest to us, and 
only the ``heftiest'' monomials
in this hierarchy matter.  This reduces the whole
problem to a combinatorial exercise.

We take $n \geq 7$ odd.
Let $m=(n-1)/2$.  Recall that
$\cal M$ is spanned by
$$O_1,...,O_m,O_n,E_1,...,E_m,E_n.$$  We define
\begin{equation}
A_{k,\pm} = O_k \pm E_k.
\end{equation}
For the first calculation, we use the monodromy invariants
\begin{equation}
\label{inv1}
A_{3,+},...,A_{m,+},A_{n,+},A_{2,-},...,A_{m,-},A_{n,-}.
\end{equation}
For the second calculation, we use the monodromy invariants
\begin{equation}
\label{inv2}
A_{3,-},...,A_{m,-},A_{3,+},...,A_{m,+},A_{n,+}.
\end{equation}

The point we use is of the form $p=P^u$, where
$P^u$ is an $n$-gon having
corner invariants
\begin{equation}
a,b,c,d,u^1,u^2,u^3,u^4,...,u^4,u^3,u^2,u^1,d,c,b,a,
\end{equation}
Here \begin{itemize}
\item $a = O(u^{(n-4)(n-3)/2})$.
\item $b = 1 +O(u)$
\item $c = 1 +O(u)$.
\item $d = 1+ O(u)$.
\end{itemize}
We will show that the results hold when $u$ is
sufficiently small.
Here we are using the big O notation, so that
$O(u)$ represents an expression that is at most $Cu$
in size, for a constant $C$ that does not depend on $u$.

We will construct $P^u$ in the next section.
Our first calculation requires only the information
presented above.  The second calculation, which
is almost exactly the same as the first calculation,
requires some auxilliary justification.  In order
to justify the calculation we make, we need to
make some estimates on the tangent space $T_{P^u}$
to $\cC_n$ at $P^u$.
We will also do this in the next section.  

In Section \ref{calc1} and Section \ref{calc2} we will explain our
two calculations in general terms.
In Section \ref{heft} we will define the concept of
the {\it heft\/} of a monomial, and we will use
this concept to put a kind of ordering on the
monomials that appear in the monodromy invariants
of interest to us.  Following the analysis of
the heft, we complete the details of our calculations.

\subsection{The first calculation in broad terms}
\label{calc1}

Let $\nabla$ denote the gradient on $\R^{2n}$.
Let $\widetilde \nabla$ denote the {\it normalized gradient\/}:
\begin{equation}
\widetilde \nabla F= \lambda^{-1}\nabla F; \hskip 30 pt
\lambda=\|\nabla F\|_{\infty}.
\end{equation}
In practice, we never end up dividing by zero.
So, the largest entry in $\widetilde \nabla F$ is
$\pm 1$.

If $F$ is a monodromy invariant, the coordinates
of $\widetilde \nabla F(P^u)$ have a power series in $u$.
We define $\Psi F$ to be the result of setting all
terms except the constant term to $0$. We call
$\Psi F$ the {\it asymptotic gradient\/}.  Thus,
if $$\widetilde \nabla F(P^u)=(1-u^3\cdots,-1+u\cdots,u^2\cdots,...)$$ then
$\Psi F = (1,-1,0,...)$.

\begin{lemma}
\label{degen}
Suppose that 
$\Psi F_1,...,\Psi F_k$ are
linearly independent.
Then likewise
$\nabla F_1,...,\nabla F_k$ are
linearly independent at $P^u$ for
$u$ sufficiently small.
Equivalently, the same goes for
$dF_1,...,dF_k$.
\end{lemma}

\startproof
Since $\Psi F_1,...,\Psi F_k$ are independent there is
some $\epsilon>0$ such that a sum of the form
$$\bigg|\sum b_j \Psi F_j\bigg|<\epsilon; \hskip 30 pt
\max |b_j|=1$$
is impossible.

Suppose for the sake of contradiction that
the gradients are linearly dependent at $P^u$
for all sufficiently small $u$.  Then the
normalized gradients are also linearly dependent
at $P^u$ for all sufficiently small $u$.  We may write
\begin{equation}
\sum b_j \widetilde \nabla F_j \cdot e_i=0;
\hskip 30 pt \max |b_j|=1.
\end{equation}
for the standard basis vectors $e_1,...,e_{2n}$.
The coefficients $b_j$ possibly depend on $u$,
but this doesn't bother us.

We have the bound
\begin{equation}
\bigg|b_j \widetilde \nabla F_j - b_j \Psi F_j\bigg|= O(u).
\end{equation}
Hence
\begin{equation}
\sum_j b_j \Psi F_j \cdot e_i=O(u)
\end{equation}
for all basis vectors $e_i$.  Therefore,
we can take $u$ small enough so that
$$\bigg|\sum b_j \Psi F_j\bigg|<\epsilon; \hskip 30 pt
\max |b_j|=1,$$
in contradiction to what we said at the beginning of
the proof.
\proofend

\begin{remark}
{\rm  The idea of the proof of the previous lemma is simple: given a matrix, algebraically dependent on a parameter $u$, the rank of the matrix is greatest in a Zariski open subset of the parameter space and can only drop for special values of the parameter (zero, in our case).}
\end{remark}

We form a matrix $M_+$ whose rows
are $\Psi F$, where $F$ is each of
the $A_+$ invariants.
We similarly form the matrix $M_-$.

\begin{lemma}
\label{ortho}
Each row of $M_+$ is orthogonal to each row of $M_-$.
\end{lemma}

\startproof
Consider the map
$T: \R^{2n} \to \R^{2n}$ which simply reverses
the coordinates.  We have $E_k \circ T=O_k$ for all $k$
and moreover $T(P^u)=P^u$.
Letting $dT$ be the differential of $T$, we have
\begin{equation}
dT(\nabla A_{k,\pm})= \pm \nabla A_{k,\pm}.
\end{equation}
Our lemma follows immediately from this equation,
and from the fact that $T$ is an isometric involution.
\proofend

In view of Lemmas \ref{degen} and Lemma \ref{ortho}, our
first calculation follows from the statements that
$M_+$ and $M_-$ have full rank.  

For the matrix $M_+$, we consider the minor
$m_+$ consisting of columns
$$1,6,7,10,11,14,15,18,19,...$$
until we have a square matrix.
We will prove below that $m_+$ has the following form
(shown in the case $n=13$.)
\begin{equation}
\label{special}
\left[\begin{matrix}
0&\pm 1&\pm 1&\pm 1&\pm 1\cr
0&0&\pm 1&\pm 1&\pm 1\cr
0&0&0&\pm 1&\pm 1\cr
0&0&0&0&\pm 1\cr
\pm 1&0&0&0&0 \end{matrix}\right]
\end{equation}
This matrix always has full rank.
Hence $M_+$ has full rank.

For the matrix $M_-$ we consider the minor
$m_-$ consisting of columns
$$1,{\bf 3\/},6,7,10,11,14,15,18,19,...$$
The only difference here is that column $3$ is
inserted.  The resulting matrix has
exactly the same structure as just
described.
Hence $M_-$ has full rank.

\subsection{The second calculation in broad terms}
\label{calc2}

Let $T=T_{P^u}(\cC_n)$ denote the tangent space
to $\cC_n$ at $P^u$.  Let
$\{e_k\}$ denote the standard basis for
$\R^{2n}$. Let $\pi: \R^{2n} \to \R^{2n-8}$ denote the map
which strips off the first and last $4$ coordinates.

Define
\begin{equation}
\nabla_8=\pi \circ \nabla.
\end{equation}
We define the normalized version 
$\widetilde \nabla_8$ exactly as we defined
$\widetilde \nabla$. Likewise we define
$\Psi_8 G$ for any monodromy function $G$.

For a collection of vectors $v_5,\dots,v_{2n-4}$ to be specified in the next lemma, we form the vector
\begin{equation}
\Upsilon_8 G=(D_{v_5}G,...,D_{v_{2n-4}}G)
\end{equation}
made from the directional derivatives of
$G$ along these vectors.  Note, by way of
analogy, that
\begin{equation}
\nabla_8 G=(D_{e_5}G,...,D_{e_{2n-4}}G).
\end{equation}
We define
the normalized version $\widetilde \Upsilon_8$ exactly
as we defined $\widetilde \nabla_8$.

In the next section, we will
establish the following result.

\begin{lemma}[Justification]
\label{just}
There is a basic
$v_5,...,v_{2n-4}$ for $T_{P^u}(\cC_n)$ 
such that
$\pi(v_k)=e_k$ for all $k$ and
$$\widetilde \Upsilon_8 G - \widetilde \nabla_8 G = O(u).
$$
\end{lemma}

\begin{corollary}
Suppose that 
$\Psi_8 G_1,...,\Psi_8 G_k$ are
linearly independent. Then
the restrictions of
$dG_1,...,dG_k$ to $T_{P^u}(\cC_n)$ are
linearly independent for $u$ sufficiently small.
\end{corollary}

\startproof
Given our basis, $\Psi_8$ represents the constant
term approximation of both
$\widetilde \Upsilon_8$
and $\widetilde \nabla_8$. So, the same proof
as in Lemma \ref{degen} shows that
the vectors
$\widetilde \Upsilon_8 G_j$ are linearly independent.
This is equivalent to the conclusion of
our corollary.
\proofend

Using the invariants listed in (\ref{inv2}),
we form the matrices $M_+$ and $M_-$ just as
above, using $\Psi_8$ in place of $\Psi$.
Lemma \ref{ortho} again shows that
each row of $M_+$ is orthogonal to
each row of $M_-$.
Hence, we can finish the second calculation by
showing that both $M_+$ and $M_-$ have full rank.

For $M_-$ we create a square minor $m_-$ using the
columns
$$2,3,6,7,10,11,14,15,...$$
Again, we continue until we have a square.
It turns out that $m_-$ has the form
\begin{equation}
\label{special2}
\left[\begin{matrix}
\pm 1&\pm 1&\pm 1&\pm 1\cr
0&\pm 1&\pm 1&\pm 1\cr
0&0&\pm 1&\pm 1\cr
0&0&0&\pm 1\cr \end{matrix}\right]
\end{equation}
Hence $M_-$ has full rank.

For $M_+$ we create a square minor $m_+$ using the
same columns, but extending out one further
(on account of the larger matrix size.)
It turns out that $m_+$ has the same form as $m_-$.
Hence $M_+$ has full rank.

\subsection{The heft}
\label{heft}

Any monomial in
the variables $x_1,...,x_{2n}$, when
evaluated at $P^u$, has a power
series expansion in $u$.  We define the
{\it heft\/} of the monomial to be the
smallest exponent that appears in this series.
For instance, the heft of $u^2+u^3$ is $2$.
We define the heft of a polynomial to be the
minimum heft of the monomials that comprise it.
Given a polynomial $F$, we define
heft of $\nabla F$ to be the minimum
heft, taken over all partial
derivatives $\partial F/\partial x_j$.

We call a monomial term of
$\partial F/\partial x_k$ {\it hefty\/} if its heft realizes
the heft of $\nabla F$.
We define $H_kF$ to be the sum of the
hefty monomials in $\partial F/\partial x_k$.
Each monomial occurs with sign $\pm 1$. We
define $|H_kF| \in \Z$ to be the sum of the
coefficients of the hefty terms in $H_kF$.
We say that $F$ is {\it good\/} if
$|H_kF| \not = 0$ for at least one index $k$.
If $F$ is good then
\begin{equation}
\label{good}
\Psi F= C (|H_1F|,...,|H_{2n}F|),
\end{equation}
for some nonzero constant $C$ that depends on $F$.
It turns out that $C=\pm 1$ in all cases.

We say that $F$ is {\it great\/} if $F$ is
{\it good\/} and $|H_k F| \not = 0$ for
at least one index $k$ which is not amongst
the first or last $4$ indices.  When $F$ 
is great, not only does equation (\ref{good})
hold, but we also have
\begin{equation}
\label{great}
\Psi_8 F= C (|H_5F|,...,|H_{2n-4}F|),
\end{equation}

\begin{lemma}
Let $k=2,3$. Then
$A_{k,\pm}$ is great and
$\nabla A_{k,\pm}$ has heft $0$.
\end{lemma}

\startproof
Let $F=A_{k,\pm}$.
Consider the case $k=2$.
The argument turns out to be the same in the
$(+)$ and $(-)$ cases. We say that an {\it outer variable\/} is one of the
first or last $4$ variables in $\R^{2n}$, and we call
the remaining variables {\it inner\/}.
Since $x_2x_6$ and $x_6x_{2n-2}$ are both terms of
$F$, we see that $$H_6F=x_2+x_{2n-2}+...$$
In particular, $\nabla F$ has heft $0$.  Any term
in $H_6F$ involves only the outer $8$ variables,
and a short case-by-case analysis shows that there
are no other possibilities besides the two terms
listed above.  Hence $|H_6F|=2$.  This shows that
$F$ is great.

Now consider the case $k=3$.
The argument turns out to be the same in the
$(+)$ and $(-)$ cases.
Since $x_2x_6x_{2n-2}$ is a term of $F$
 we see that $$H_6F=x_2x_{2n-2}+...$$
The rest of the proof is as in the previous
case, with the only difference being
that $|H_6F|=1$ in this case.
\proofend

From now on, we fix some $F=A_{k,\pm}$ with $3<k \leq m$.
Let $\alpha_1,\alpha_2,...$ be the terms of
the following sequence
\begin{equation}
0,0,0,2,3,6,7,10,11,14,15,...
\end{equation}

\begin{lemma}
\label{array}
$\nabla F$ has heft 
at most $\alpha_1+...+\alpha_k$.
\end{lemma}

\startproof
We describe a specific term in $\nabla F$ having heft
$\alpha_1+...+\alpha_k$.  We make a monomial
using the indices
\begin{equation}
\label{packing}
2,2n-2,6,2n-6,10,2n-10,...
\end{equation}
stopping when we have used $k-1$ numbers.
The monomial corresponding to these
indices has heft
$$0+0+0+2+3+6+7+10+11...=\alpha_1+...+\alpha_k.$$
Thinking of our indices cyclically, we see that
our integers lie in an interval of length
$4k-7$. So, between the largest index in
(\ref{packing}) that is less than
$n$ and the smallest index greater than $n$ there
is an unoccupied stretch of at least $9$ integers.
The point here is that
$$9 + (4k-7) \leq 9 + 4m-7 = 9 + 2(n-1)-7=2n.$$

Given that the unoccupied stretch has at least
$9$ consecutive integers, there is at least
$1$ (and in fact at least $2$) even indices
$j$ such that the monomial
$$m=\pm x_j x_2 x_{2n-2} x_6 x_{2n-6}x_{10}...$$
is a term of $F$.
But then $\partial m/\partial x_j$ has heft
$\alpha_1+...+\alpha_k$.
\proofend

We mention that  (\ref{packing}) is one
of two obvious ways to make a term of
heft $\alpha_1+...+\alpha_k$. The other
way is to take the {\it mirror image\/}, namely
\begin{equation}
\label{packing2}
2n-1,3,2n-5,7,2n-9,11,...
\end{equation}

\begin{lemma}
\label{inner}
If $\partial F/\partial x_j$ has a hefty term, then
$j$ is an inner variable.
\end{lemma}

\startproof
For ease of exposition, we will consider the
case when $j$ is one of the first $4$ variables.
Let $(i_1,...,i_d)$ be the sequence of indices which
appear in a term $m'$ of $\partial F/\partial x_j$.
The corresponding term $m$ in $F$ has
index sequence $(j,i_1,...,i_d)$, where these
numbers are not necessarily written in order.
We know that at least one of the indices, say $a$,
is an inner variable.  By construction
$\partial m/\partial x_a$ has smaller heft than
$m'$.  Hence $\partial F/\partial x_j$ has
no hefty terms.  Hence $j$ is an inner variable.
\proofend

\begin{lemma}
Suppose 
the monomial $\pm x_{i_1}...x_{i_a}$ is a
hefty term of $\partial F/\partial x_j$.
Then $a=k-1$ and
$i_1,...,i_{k-1}$ are either
as in  (\ref{packing}) or as in
equation (\ref{packing2}).
\end{lemma}

\startproof
We have to play the following game: We have
a grid of $2n$ dots.  The first and last dot
are labelled $(n-3)(n-4)/2$.  The remaining
$6$ outer dots are labelled $0$.  The inner
dots are labelled $1,2,3,...,3,2,1$.
Say that a {\it block\/} is a collection
of $d$ dots in a row for $d=1,2,3$.  We
must pick out either $k$ or $k-1$ blocks
in such a way that the total sum of the
corresponding dots is as small as possible,
and the (cyclically reckoned) 
spacing between consecutive blocks
is at least $4$.  That is, at least $3$
``unoccupied dots'' must appear between
every two blocks.  

It is easy to see that one should use $k-1$
blocks, all having size $1$.  Moreover,
half (or half minus one) of the blocks should
crowd as much as possible to the left and
half minus one (or half) of the blocks
should crowd as much as possible to the
right. A short case by case analysis of
the placement of the first and last blocks
shows that one must have precisely the
choices made in  (\ref{packing})
and (\ref{packing2}).
\proofend

\begin{corollary} 
\label{GREAT}
Let $F=A_{k,\pm}$, with $k \geq 2$.  Then $F$ is good.
If $k \leq m$ then $F$ is great, and the heft of $\nabla F$ 
is $\alpha_1+...+\alpha_k$.
\end{corollary}

\startproof
In light of the results above, the only nontrivial
result is that $F$ is great when $3<k \leq m$.  The construction
in connection with  (\ref{packing}) produces
a hefty term of $\partial F/\partial x_j$ for
some inner index $j$.  The key observation is that,
for parity considerations, the mirror term corresponding
to  (\ref{packing2}) is not a term of
$\partial F/\partial x_j$. In one case $j$ must
be odd and in the other case $j$ must be even.
Hence, there is only $1$ hefty term in
$\partial F/\partial x_j$. 
\proofend

As regards the heft,
we have done everything but analyze the Casimirs.
Recall that
\begin{equation}
O_n=x_1x_3...x_{2n-1}; \hskip 30 pt
E_n=x_2x_4...x_{2n}.
\end{equation}

\begin{lemma}
\label{cas1}
$A_{n,\pm}$ is good and $\nabla A_{n,+}$ has heft
$$\frac{(n-3)(n-4)}{2}.$$
Moreover,
$$\Psi A_{n,\pm} = (1,0,...,0,\pm 1).$$
\end{lemma}

\startproof
Let $F$ be either of these functions.
Clearly the hefty terms of
$\nabla F$ are the ones which omit
the first and last variables.
From here, this lemma is an exercise in arithmetic.
\proofend

A similar argument proves
\begin{lemma}
\label{cas2}
$A_{n,\pm}$ is good and
$\nabla_8 A_{n,+}$ 
has heft $(n-4)^2$.  Moreover,
$$\Psi A_{n,\pm} = (0,...,0,\pm 1,1,0,...,0),$$
with the $2$ middle indices being nonzero.
\end{lemma}

\subsection{Completion of the first calculation}

To complete the first calculation, we need to analyze the
matrix made from the asymptotic gradients
$\Psi F_1,\Psi F_2,...$.  We deal with the first
two in a calculational way.

\begin{lemma}
\label{triv1}
$\Psi A_{2,\pm }=(0,0,\pm 1,0,0,1,\pm 1,...,1,\pm 1,0,0,1,0,0)$.
\end{lemma}

\startproof
Let $F=A_{2,\pm}$.  We know that $F$ has heft $0$,
so the hefty terms in $\nabla F$ are monomials
which only involve the outer indices.  Hence,
when $8 \leq j \leq 2n-8$ the result only depends
on the parity of $j$ and neither the value of $j$
nor the value of $n$.  For the remaining indices,
the result is also independent of $n$. Thus, a
calculation in the case (say) $n=13$ is general
enough to rigorously establish the whole pattern.
This is what we did.
\proofend

\begin{lemma}
\label{triv2}
$\Psi A_{3,\pm }=(0,0,0,0,0,1,\pm 1,...,1,\pm 1,0,0,1,0,0)$.
\end{lemma}

\startproof
Same method as the previous result.
\proofend

Now we are ready to analyze the minors $m_+$ and $m_-$
described in connection with the first calculation.
When we say that a certain part of one of these
matrices has the form given by  (\ref{special}),
we understand that  (\ref{special}) gives a
smallish member of an infinite family of matrices,
all having the same general type. So, we mean to
take the corresponding member of this family which
has the correct size.  

We say that a given
row or column of one of our matrices {\it checks\/} if
it matches the form given by  (\ref{special}).
We will give the argument for $m_+$.  The case for $m_-$ is
essentially the same.

\begin{lemma}
The first column of $m_+$ checks.
\end{lemma}

\startproof
By Lemma \ref{cas1},
the first coordinate of $\Psi A_{n,+}$ is $\pm 1$.
By Lemmas \ref{inner}, \ref{triv1},
and \ref{triv2}., we have
$\Psi A_{k,+}$ is zero for $k<n$. This is 
equivalent to the lemma.
\proofend

\begin{lemma}
The first row of $m_+$ checks and the last row of
$m_+$ checks.
\end{lemma}

\startproof
The first statement follows immediately from Lemma \ref{triv2}.
The second statement follows immediately from
Lemma \ref{cas1}.
\proofend

Now we finish the proof.
Consider the $i$th row of $m_+$. 
Let $k=i+2$. In light of the trivial
cases taken care of above, we can assume that
$3<k \leq m$.
Let $F=A_{k,+}$. As we discussed
in the proof of Corollary \ref{GREAT}, each
polynomial $\partial A/\partial F_j$ has
either $0$ or $1$ hefty terms. 

Assume that $j$ is even.  Let $J \subset \{1,...,2n\}$
be the unoccupied stretch from Lemma \ref{array}.
Let $J' \subset J$ denote the smaller set obtained
by removing the first and last $3$ members
from $J$.  It follows from the construction
in Lemma \ref{array} that $\partial F/\partial j$
has a hefty term if and only if $j \in J'$.
Thus the $j$th entry of the $k$th row is
$\pm 1$ if and only if $j \in J'$.  Similar
considerations hold when $j$ is odd.  
It is an exercise to show that the conditions
we have given translate precisely into
the form given in  (\ref{special}).
Hence $m_+$ checks.

\begin{remark} 
{\rm One can approach the proof differently.
When we move from row $k$ to row $k+2$ the
corresponding interval $J'=(a,b)$ changes to the
new interval $J'=(a+4,b-4)$.  From this fact, and
from our choice of minors, it
follows easily that row $k$ checks if and only
if row $k+2$ checks.  At the same time, when $n$
is replaced by $n+2$, the interval $J'=(a,b)$ changes
to $J'=(a,b+4)$.  This translates into the statement
that row $k$ checks for $n$ if and only if
row $k$ checks for $n+2$.  All this reduces the
whole problem to a computer calculation of the
first few cases. We did the calculation up to the
case $n=13$ and this suffices.}
\end{remark}

\subsection{Completion of the second calculation}

We make all the same definitions and conventions for
the second calculation, using the matrix (family) in
 (\ref{special2}) in place of the
matrix (family) in (\ref{special}).
The argument for the second calculation is really
just the same as the argument for the first
calculation.  Essentially, we just ignore the
outer $8$ coordinates and see what we get.
What makes this work is that all the functions
except $A_{n,\pm}$ are great -- the inner
indices determine the heft.  To handle
the last row of $m_+$, which involves the
Casimir $A_{n,+}$, 
we use Lemma \ref{cas2} in place of
Lemma \ref{cas1}. 

It remains to establish the Justification Lemma \ref{just}.
It is convenient to define
\begin{equation}
\delta=\frac{(n-4)(n-5)}{2}.
\end{equation}
We also mention several other pieces of notation
and terminology.  
When we line up the indices
$5,...,2n-4$, there are $2$ {\it middle indices\/}.
When $n=7$ the middle indices of $5,6,7,8,9,10$ are
$7$ and $8$. 
Let $\pi^{\perp}$ denote the projection
from $\R^{2n}$ onto $\R^8$ obtained by stringing
out the first and last $4$ coordinates.

\begin{lemma}[Tangent Estimate]
\label{est}
The following properties of $\pi^{\perp}(v_j)$ hold: 
\begin{itemize}
\item All coordinates are $O(1)$.
\item Coordinates $3$ and $6$ are $O(u)$.
\item Except when $j$ is one of the middle two indices,
coordinates $1$ and $8$ are $O(u^{\delta+1})$.
\item When $j$ is the first middle index,
coordinate $1$ is $u^{\delta}+O(u^{\delta+1})$ and
coordinate $8$ is $O(u^{\delta+1})$.
\item When $j$ is the second middle index,
coordinate $8$ is $u^{\delta}+O(u^{\delta+1})$ and
coordinate $1$ is $O(u^{\delta+1})$.
\end{itemize}
\end{lemma}

\startproof
We prove this in the next section.
\proofend

\begin{lemma}
The Justification Lemma holds for $F=A_{n,+}$.
\end{lemma}

\startproof
A direct calculation shows that, up to
$O(u^{\delta+1})$,
\begin{equation}
\widetilde \nabla F = 
(1,0,...,0,u^{\delta},u^{\delta},0,...,0,1)
\end{equation}
Hence
\begin{equation}
\widetilde \nabla_8 F=(0,...,0,1,1,0,....0)+ O(u).
\end{equation}

Let $Z$ be the first coordinate of $\nabla F$.
If $j$ is not a middle index, we have
\begin{equation}
D_{v_j}F=\nabla F \cdot v_j=Z \times O(u^{\delta+1}).
\end{equation}
This estimate comes from the Tangent Estimate Lemma \ref{est}.

If $j$ is the first middle index, then
\begin{equation}
D_{v_j}F=\nabla F \cdot v_j=Z \times 2 O(\delta).
\end{equation}
The first contribution comes from coordinate 1,
and is justified by the Tangent Estimate Lemma,
and the second contribution comes from coordinate $j$.

The above calculations show that
\begin{equation}
\widetilde \Upsilon_8 F= (0,...,0,1,1,0,....0)+ O(u).
\end{equation}
Hence $\widetilde \nabla_8 F = \widetilde \Upsilon_8 F+O(u)$.
\proofend

Now suppose that $F$ is one of the relevant monodromy
invariants, but not the Casimir.

Our analysis  establishes

\begin{lemma}
\label{grad}
Both $\pi^{\perp}(\widetilde \nabla F)$ and
$\pi^{\perp}(\nabla F)$ have the following properties.
\begin{enumerate}
\item All coordinates 
 are at most $1+O(u)$ in size.
\item All coordinates 
 except coordinates $3$ and $6$ are $O(u)$.
\end{enumerate}
\end{lemma}

\startproof
This is immediate from our analysis of the heft
of $\nabla F$.
\proofend

\begin{lemma} One has
$$\widetilde \nabla_8 F \cdot e_j =
\widetilde \nabla F \cdot v_j + O(u).
$$
\end{lemma}

\startproof
Combining the Tangent Estimate Lemma with Lemma \ref{grad},
we see that
$$
\pi^{\perp}(\widetilde \nabla F) \cdot \pi^{\perp}(v_j)=O(u).
$$
Hence
\begin{equation}
\label{zoop1}
\widetilde \nabla F \cdot v_j = \pi \circ \widetilde \nabla F
\cdot e_j  + O(u).
\end{equation}

From Property 1 above, we see that
$$\|\nabla_8F\|_{\infty}=
\|\nabla F\|_{\infty}+O(u).$$
Therefore
\begin{equation}
\label{zoop2}
\widetilde \nabla_8 F=\pi \circ \widetilde \nabla F+ O(u).
\end{equation}

Combining equations (\ref{zoop1}) and (\ref{zoop2}), we get
the result of the lemma.
\proofend

\begin{lemma}
Setting
$\lambda=\|\nabla F\|_{\infty}$, we have
$$
\lambda^{-1}(\Upsilon_8 F)_j = (\widetilde \Upsilon_8 F)_j + O(u).
$$
Here $(X)_j$ is the $j$th coordinate of $X$.
\end{lemma}

\startproof
Combining the Tangent Estimate Lemma \ref{est} with Lemma \ref{grad}, we have
$$\pi^{\perp} \circ \nabla F \cdot \pi^{\perp}(v_j) = O(u).$$
Therefore
$$
\|\Upsilon_8 F\|_{\infty} = \|\nabla_8 F\|_{\infty}+O(u).
$$
Combining this with equation (\ref{zoop2}), we have
$$
\|\Upsilon_8 F\|_{\infty} =\|\nabla F\|_{\infty} + O(u).
$$
Our lemma follows immediately.
\proofend

By definition, we have
\begin{equation}
\widetilde \nabla F \cdot v_j =
\lambda^{-1} \nabla F \cdot v_j=\lambda^{-1}(\Upsilon_8 F)_j; \hskip 30 pt
\lambda=\|\nabla F\|_{\infty}.
\end{equation}
Combining this last equation with our two lemmas, we have
\begin{equation}
(\widetilde \nabla_8 F)_j=\widetilde \nabla_8 F \cdot e_j =
(\widetilde \Upsilon_8 F)_j + O(u).
\end{equation}
This holds for all $j$.  This completes the proof of
the Justification Lemma.

\section{The polygon and its tangent space}
\label{TAN}

The goal of this section is to construct the polygon
$P^u$ and prove the Tangent Lemma, which estimates
the tangent space $T_{P^u}(C)$.   We will begin by
repackaging some of the material worked out
in \cite{Sch3}. The results here are self-contained,
though our main formula relies on the work done
in \cite{Sch3}.   In order to remain consistent with
the formulas in \cite{Sch3}, we will use a slightly
different labelling convention for polygons.

\subsection{Polygonal rays}

We say that a {\it polygonal ray\/} is an infinite
list of points $P_{-7}, P_{-3}, P_1,P_5,...$ in the projective plane.
We normalize so that (in homogeneous coordinates)
\begin{equation}
P_{-7}=(0,0,1), \hskip 30 pt
P_{-3}=(1,0,1), \hskip 30 pt
P_1 = (1,1,1), \hskip 30 pt
P_5=(0,1,1).
\end{equation}
The first $4$ points are normalized to be the vertices of
the positive unit square, starting at the origin, and
going counterclockwise.  Here we are interpreting these
points in the usual affine patch $z=1$.
This polygonal ray defines lines:

\begin{equation}
L_{-5+k}=P_{-7+k}P_{-3+k}; \hskip 30 pt k=0,4,8,...
\end{equation}
We denote by $LL'$ the intersection $L \cap L'$.
Similarly, $PP'$ is the line containing $P$ and $P'$.

The pairs of points and lines determine flags, as follows:
\begin{equation}
F_{-6+k}=(P_{-7+k},L_{-5_k}), \hskip 30 pt
F_{-4+k}=(P_{-3+k},L_{-5+k}), \hskip 30 pt
k=0,4,8,12...
\end{equation}

The corner invariants were defined in Section \ref{CORNERDEF}.
In this section we relate the definition there to
our labelling convention here.
We define

\begin{equation}
\label{flag2}
\chi(F_{0+k})=[P_{-7+k},P_{-3+k},L_{-5+k}L_{3+k},L_{-5+k}L_{7+k}],
\hskip 30 pt k=0,4,8,...
\end{equation}

\begin{equation}
\chi(F_{2+k})=[P_{9+k},P_{5+k},L_{7+k}L_{-1+k},L_{7+k}L_{-5+k}],
\hskip 30 pt k=0,4,8,...
\end{equation}
Here we are using the inverse cross ratio, as in equation
\ref{ICR}.
Referring to the corner invariants, we have
\begin{equation}
x_k=\chi(F_{2k}); \hskip 30 pt
x_{k+1}=\chi(F_{2k+2}); \hskip 40 pt k=0,2,4,...
\end{equation}

\begin{remark}
{\rm Notice that it is impossible to define $\chi(F_{-2})$ because we
would need to know about a point $P_{-11}$, which we have
not suppled.  Likewise, it is impossible to define
$\chi(F_{-4})$ because we would need to know about $L_{-9}$,
which we have not supplied.  Thus, the invariants
$x_0,x_1,x_2,...$ are well defined for our polygonal ray.}
\end{remark}

{\bf  Cross product in vector form:\/}
Since we are going to be computing a lot of these
cross ratios, we mention a formula that works
quite well.  We represent both points and lines
in homogeneous coordinates, so that
$(a,b,c)$ represents the line corresponding
to the equation $ax+bx+cz=0$.  We define
$V*W$ to be the coordinate-wise product of
$V$ and $W$.  Of course, $V*W$ is also a
vector.  Let $(\times)$ denote the cross
product.  We have
\begin{equation}
\label{cross}
(\chi,\chi,\chi)=\frac{(A \times B)*(C \times D)}{(A \times C)*(B \times D)}.
\end{equation}
Here $\chi$ is the inverse cross ratio of the points or
lines represented by these vectors.
It may happen that some coordinates in the denominator vanish.
In this case, one needs to interpret this equation as a kind
of limit of nearby perturbations.
This formula works whenever $A,B,C,D$ represent either
collinear points or concurrent lines in the projective plane.

\subsection{The reconstruction formulas}
\label{RECON}

Referring to the definition of the monodromy invariants, we
define $O_a^b$ to be the sum over all odd admissible monomials
in the variables $x_0,x_1,x_2,...$ which do not involve
any variables with indices $i \leq a$ or $i \geq b$.    For instance
$$O_1^1=1, \hskip 20 pt O_1^3=1, \hskip 20 pt
O_1^5=1-x_3, \hskip 20 pt O_1^7 =  1 -x_3 +x_3x_4x_5.$$
We also note that, when $a<0$, the polynomial
$O_a^b$ is independent of the value of $a$.
For this reason,
when $a<0$ we simply write $O^b$ in place of $O_a^b$.
The corresponding set $S^b$ consists of admissible
sequences, all of terms are less than $b$.

Given a list $(x_0,x_1,x_2,...)$ we seek a polygonal ray which
has this list as its corner invariants.  Here is the formula.

\begin{equation}
\label{reconstruct}
P_{9+2k}=(O^{3+k} - O_1^{3+k} + x_0x_1 O_3^{3+k},O^{3+k},O^{3+k} + x_0x_1 O_3^{3+k}),
\hskip 15 pt k=0,2,4,...
\end{equation}

We would also like a formula for reconstructing the lines
of a polygonal ray. 
We start with the obvious:
\begin{equation}
L_{-5}=(0,1,0); \hskip 30 pt
L_{-1}=(-1,0,1); \hskip 30 pt
L_3=(0,-1,1).
\end{equation}
For the remaining points, we define
polynomials
$E_a^b$ exactly as we defined $O_a^b$ except we 
interchange the uses of {\it even\/} and {\it odd\/}.
Thus, for instance $E_2^6=1-x^4$.  Here is the formula.
\begin{equation}
\label{reconstruct3}
L_{7+2k} = (E^{2+k}-E_0^{2+k},E_0^{2+k}-x_0 E_2^{2+k},-E^{2+k}),
\hskip 30 pt k=0,2,4...
\end{equation}

\begin{remark}
{\rm These formulas are equivalent to equations 19
and 20 in \cite{Sch3}, but the normalization
of the first $4$ points is different, and
the roles of points and lines have been switched.
We got the above formulas by applying a suitable
projective duality to the polygonal ray in
\cite{Sch3}. }
\end{remark}

We mention one important connection between our
various reconstruction formulas.  The
following is an immediate consequence of
Lemma 3.2 in \cite{Sch3}:
\begin{equation}
\label{reconstruct4}
P_{5+k} \times P_{9+k} = -(x_1x_3x_5,...,x_{k/2+1}) L_{7+k},
\hskip 20 pt k=0,4,8...
\end{equation}

We close this section with a characterization of the 
moduli space of closed polygons within $X$.  We
do not need this result for our proofs, but it
is nice to know.\footnote{One could give an alternative proof of Proposition \ref{TanProp} computing the Poisson bracket of the polynomials of Lemma \ref{variety} with the monodromy invariants.} 
\begin{lemma}
\label{variety}
The invariant $x_1,...,x_{2n}$ define a closed polygon
if and only if $O^{2n-5}$ and all its cyclic shifts
vanish.
\end{lemma}

\startproof
We can think of a closed polygon as an $n$-periodic
infinite ray.  The periodicity implies that
$P_{4n-7}=P_{-7}=(0,0,1)$.  Since
$4n-7 = 2k+9$ for $k=2n-8$, 
equation (\ref{reconstruct}) tells us that
$O^{2n-5}=0$.  Considering equation
\ref{reconstruct3}, we see that
$E^{2n-6}=E_0^{2n-6}=0$.  But $E_0^{2n-6}$ is a cyclic
shift of $O^{2n-5}$.  Hence, if $P$ is closed then
$O^{2n-5}$ and all its cyclic shifts vanish.

Conversely, if $O^{2n-5}$ and all its shifts vanish then
$P_{4n-7} \in L_{-5}$ and $P_{-3} \in L_{4n-5}$.  Likewise
$P_{4n-3} \in L_{-1}$ and $P_{1} \in L_{4n-1}$, and so on.
This situation forces $P_{4n-3}=P_{-3}$.  Shifting
the indices, we see that $P_{4n+1}=P_1$, and so on.
\proofend

\begin{remark}
{\rm Observe that $O^{2n-5}$ involves exactly $2n-7$
consecutive corner invariants.  If the
first $2n-8$ are specified, then the next
variable can be found by solving $O^{2n-5}=0$.
Thus, Lemma \ref{variety} gives an algorithmic
way to find a closed $n$-gon whose first
$2n-8$ corner invariants are specified.}
\end{remark}

\subsection{The polygon}

We start with an infinite periodic list of variables
which starts out
\begin{equation}
\label{corner}
(u,u^2,...,u^{n-4},u^{n-4},...,u^2,u^1,...)
\end{equation}
and has period $2n-8$.
We let $X_u$ denote the polygonal ray associated to this
infinite list.  Once $u$ is sufficiently small, the
first $n$ points of $X_u$ are well defined.  We define $P^u$ to be the $n$-gon
made from the first $n$-points of $X_u$, and we take
$u$ small enough so that this definition makes sense.

The first $2n-8$ corner invariants of $P^u$,
which we now identify with $x_0,...,x_{2n-9}$, are
the ones listed in equation (\ref{corner}).
However, when it comes time to compute
$x_{2n-8},...,x_{2n-1}$, we do not use the
relevant points of $X_u$ but rather substitute
in the corresponding point of $P^u$.
Thus, the remaining $8$ corner invariants change.
We write the corner invariants of $P^u$ as
\begin{equation}
a,b,c,d,u,u^2,u^3,...,u^3,u^2,u,d',c',b',a'.
\end{equation}

It follows 
from symmetry that $e=e'$ for each $e \in \{a,b,c,d\}$.
This symmetry here is that the
first $2n-8$ invariants determine $P$, and their
palindromic nature forces $P$ to be
self-dual: the projective duality carries $P$
to the dual polygon made from the lines
extending the sides of $P$. 

\begin{lemma}
\label{ccc}
$e=1+O(u)$ for each $e \in \{b,c,d\}$.
\end{lemma}

\startproof
We set $P_{-11}=(X,Y,Z)$ and $L_{-13}=(U,V,W)$.  We have
\begin{equation}
L_{-9}=(1,0,0) \times (X,Y,Z)=(-Y,X,Z).
\end{equation}
Equations \ref{reconstruct} and \ref{reconstruct3} tell us
\begin{equation}
\label{orders}
(X,Y,Z)=(1,0,0)+O(u); \hskip 30 pt (U,V,W)=(0,1,-1)+O(u).
\end{equation}
We compute 
\begin{equation}
\label{go3}
b=\chi(F_{-6})=\chi(P_{1},P_{-3},L_{-1}L_{-9},L_{-1}L_{-13})=
\frac{UX+WX+VY}{(U+W)(X-Y)}.
\end{equation}

\begin{equation}
\label{go2}
c=\chi(F_{-4})=\chi(P_{-11},P_{-7},L_{-9}L_{-1},L_{-9}L_3)=
 \frac{X-Y}{X-Z}
\end{equation}

\begin{equation}
\label{go1}
d=\chi(F_{-2})=\chi(P_{5},P_{1},L_{3}L_{-5},L_{3}L_{-9})=
d=\frac{X}{X+Y+Z}.
\end{equation}

Our result is immediate from these formulas
and from equation (\ref{orders}).
\proofend

\begin{lemma}
\label{prefinal}
$a=u^s + O(u^{s+1})$, where $s=(n-4)(n-3)/2$.
\end{lemma}

\startproof
We have
\begin{equation}
\label{go4}
a=\chi(F_{-8})=\chi(P_{-15},P_{-11},L_{-13}L_{-5},L_{-13}L_{-1})=
\chi(A,B,C,D).
\end{equation}
We will estimate $a$ by considering the middle coordinate of
equation (\ref{cross}).
Calculations similar to the ones above give
$$
A = (0,1,1) + O(u), \hskip 10 pt
B = (0,1,1) + O(u), $$
\begin{equation}
C = (1,0,0) + O(u), \hskip 10 pt
D = (1,1,1) + O(u).
\end{equation}
Hence
\begin{equation}
\label{prefinal1}
(A \times C)_2 = + 1+O(u); \hskip 10 pt
(B \times D)_2 = + 1+O(u); \hskip 10 pt
(C \times D)_2 = - 1+O(u).
\end{equation}
Recall that
\begin{equation}
P_{-15}=P_{-15+4n}; \hskip 30 pt
P_{-11}=P_{-11+4n}.
\end{equation}
According to equation (\ref{reconstruct4}), we have
$$A \times B = -(x_1x_3,...x_{2n-9}) L_{-13+2n} = $$
\begin{equation}
\label{prefinal2}
-u_2 u_4 ... u_3 u_1 L_{-13+2n} =
-u^s L_{-13+4n}.
\end{equation}
But
\begin{equation}
L_{-13+4n} = (0,1,-1) + O(u).
\end{equation}
Therefore
$$(A \times B)_2=-u^s + O(u^{s+1}).$$
Looking at the signs in equation (\ref{prefinal1}), we see that
$a=u^s + O(u^{s+1})$.
\proofend

\subsection{The tangent space}

 Recall that $\pi: \R^{2n} \to \R^{2n-8}$
is the projection which strips off the outer $4$ coordinates.
Let $\pi^{\perp}$ be as in the Tangent Estimate Lemma \ref{est}.
Recall that $\{v_k\}$ is the special basis of $T_{P}(C)$
such that $\pi(v_k)=e_k$ for $k=5,...,2n-4$.

\begin{lemma}
The following holds concerning the coordinates
of $\pi^{\perp}(v_j)$:
\begin{itemize}
\item Coordinates $2,4,5,7$ of $\pi^{\perp}(v_j)$ have
size $O(1)$ 
\item Coordinates $3,6$ have size
$O(u)$.
\end{itemize}
\end{lemma}

\startproof
As above, we will just consider coordinates $2,3,4$.
The other cases follow from symmetry.

We refer to the quantities used in the proof of Lemma \ref{ccc}.
Each of these quantities is a polynomial in the coordinates,
depending only on $n$.  Hence $dX/dt$, etc., are all of size
at most $O(1)$.  Moreover, the denominators on the right
hand sides of equations (\ref{go3}), (\ref{go2}), and (\ref{go1})
are all $O(1)$ in size.  Our first claim now follows from
the product and quotient rules of differentiation.

For our second claim, we differentiate equation (\ref{go2}):
$$
\frac{dc}{dt}=\frac{X'(Y-Z) - X(Y'-Z') + ZY' - YZ'}{(X+Z)^2}=^*
$$
\begin{equation}
\label{last}
Y' - Z' + O(u)=
\frac{d}{dt} (-x_0x_1) O_3^{3+k}.
\end{equation}
The starred equality comes from the fact that
$(X,Y,Z)=(0,1,1)+O(u)$.
The claim now follows from the fact
that $x_0(0)=u$ and $x_1(0)=u^2$ and
$(O_3^{3+k})'(0)=O(1)$.
\proofend

\begin{lemma}
The following holds concerning the coordinates
of  $\pi^{\perp}(v_j)$:
\begin{itemize}
\item When $j$ is not a middle index,
coordinates 1 and 8 of are of
size $O(u^{\delta+1})$.
\item  When $j$ is the first middle index,
coordinate 1 equals $u^{\delta}(1+O(u))$ and
coordinate 8 is of size  $O(u^{\delta+1})$.
\item  When $j$ is the second middle index,
coordinate 8 equals $u^{\delta}(1+O(u))$ and
coordinate 1 is of size  $O(u^{\delta+1})$.
\end{itemize}
\end{lemma}

\startproof
We will just deal with coordinate 1.
The statements about coordinate 8 follow from symmetry.

Let us revisit the proof of Lemma \ref{prefinal}.
Let $f=-(A \times B)_2$.
We have $a=fg$, where
\begin{equation}
g=-\frac{(C \times D)_2}{(A \times C)_2(B \times D)_2}.
\end{equation}
We imagine that we have taken some variation, and all these
quantities depend on $t$. 

Each of the factors in the equation for $g$ has
derivative of size $O(1)$.  Moreover, the denominator
in $g$ has size $O(1)$.  From this, we conclude that
\begin{equation}
g(0)=1+O(u); \hskip 30 pt g'(0)=O(1).
\end{equation}
It now follows from the product
rule that
\begin{equation}
\frac{da}{dt}=\frac{df}{dt}(1+O(u)).
\end{equation}

Equations \ref{prefinal1} and \ref{prefinal2} tell us that
\begin{equation}
f(t)=(x_1x_3,...,x_{2n+9})\lambda(t); \hskip 30 pt
\lambda(t)=(L_{-13+2n})_2.
\end{equation}

By equation (\ref{reconstruct3}), we have
\begin{equation}
\lambda(0)=1+O(u); \hskip 30 pt
\lambda'(0)=O(1).
\end{equation}
Hence, by the product rule,
\begin{equation}
\frac{da}{dt} = \frac{d}{dt}(x_1x_3...x_{2n+9})(1+O(u)).
\end{equation}
Using the variables
\begin{equation}
x_1=u,...,x_j=u^j+t,x_{j+1}=u^{j+1},...
\end{equation}
we get the result of this lemma as a simple exercise in
calculus.
\proofend

The results above combine to prove the Tangent Space Lemma.

\bigskip

\noindent \textbf{Acknowledgments}. 
Some of this research was carried out in
May, 2011, when all three authors were
together at Brown University.  We would
like to thank Brown for its hospitality during this period.
ST was partially supported by a Simons Foundation grant. 
RES was partially supported by N.S.F. Grant DMS-0072607,
and by the Brown University 
Chancellor's Professorship.

\vskip 1cm


Valentin Ovsienko: 
CNRS, Institut Camille Jordan, 
Universit\'e Lyon 1, Villeurbanne Cedex 69622, France, 
ovsienko@math.univ-lyon1.fr

\medskip

Richard Evan Schwartz: 
Department of Mathematics,
Brown University,
Providence, RI 02912, USA, 
res@math.brown.edu

\medskip

Serge Tabachnikov: 
Department of Mathematics,
Pennsylvania State University,
University Park, PA 16802, USA, 
tabachni@math.psu.edu

\end{document}